\documentclass[11pt]{article}
\usepackage{amssymb,amsmath,upgreek,amsthm}
\usepackage{graphicx}
\usepackage{color,centernot}
\usepackage{fullpage}

\newcommand\lcm{{\mathrm{lcm}}}

\newtheorem{theorem}{Theorem}
\newtheorem*{theorem-1-paired-CIS}{Theorem~\ref{thm:1-paired-CIS}}
\newtheorem{lemma}{Lemma}
\newtheorem{claim}{Claim}
\newtheorem{corollary}{Corollary}
\newtheorem{proposition}{Proposition}[section]

\newtheorem{definition}{Definition}

\newtheorem{remark}{Remark}
\newtheorem*{remarks}{Remarks}
\newtheorem*{examples}{Examples}
\newtheorem*{example}{Example}

\def\ZZ{\mathbb{Z}}

\title{$\mbox{On CIS Circulants}$\thanks
{This work is supported in part by the Slovenian Research Agency
(research program P$1$--$0285$ and research projects
J$1$--$4010$, J$1$--$4021$, BI-US/$12$--$13$--$029$ and N$1$--$0011$:
GReGAS, supported in part by the European Science Foundation).
The first author also thanks for partial support the National Science Foundation (Grants CMMI-0856663 and IIS-1161476).}}

\author{
Endre Boros and Vladimir Gurvich\\
\small RUTCOR, Rutgers University, New Jersey, USA\\
\small 640 Bartholomew Rd, Piscataway NJ 08854-8003, USA\\
\small \texttt{\{boros,gurvich\}@rutcor.rutgers.edu}\\
\and
Martin Milani\v c\\
\small University of Primorska, UP IAM, Muzejski trg 2, SI6000 Koper, Slovenia\\
\small University of Primorska, UP FAMNIT, Glagolja\v ska 8, SI6000 Koper, Slovenia\\
\small \texttt{martin.milanic@upr.si}
}

\date{\today}
\begin{document}
\maketitle

\begin{abstract}
A {\it circulant} is a Cayley graph over a cyclic group.
A {\em well-covered} graph is a graph in which
all maximal stable sets are of the same size  $\alpha = \alpha(G)$,
or in other words, they are all maximum.
A {\it CIS} graph is a graph in which
every maximal stable set and every maximal clique intersect.
It is not difficult to show that
a circulant  $G$  is a CIS graph if and only if
$G$  and its complement $\overline G$  are both well-covered
and the product  $\alpha(G) \alpha(\overline G)$
is equal to the number of vertices.
It is also easy to demonstrate that
both families, the circulants and the CIS graphs, are
closed with respect to the operations of taking the complement
and lexicographic product.
We study the structure of the CIS circulants.
It is well-known that all  $P_4$-free graphs are CIS.
In this paper, in addition to the simple family
of the $P_4$-free circulants, we construct
a non-trivial sparse but infinite family of CIS circulants.
We are not aware of any CIS circulant that
could not be obtained from graphs in this family
by the operations of taking the complement and lexicographic product.

\medskip

{\bf Keywords:} Circulant, CIS graph, well-covered graph,
maximal stable set; maximum stable set; maximal clique; maximum clique.

\small {\bf Math.~Subj.~Class.}~(2010): 05C25, 05C69
\end{abstract}

\section{Introduction}\label{sec:intro}

\subsection{Basic concepts and operations}
We consider finite non-directed graphs without loops and multiple edges.
A graph $G = (V,E)$  has vertex-set $V$  and edge-set  $E$;
furthermore, $n = |V|$  and  $m = |E|$  are called
the {\em order} and {\em size} of $G$, respectively.
The {\it complement} ${\overline{G}}$  of a graph $G = (V,E)$ is
the graph with the same vertex-set  $V$  and the complementary edge-set
$\overline E = \{\{x,y\}\,\mid \,x,y\in V,~x\neq~y, \textrm{ and } \{x,y\}\not\in E\}$.
We say that a graph $G$ is {\it co-connected} if its complement is connected.
A {\em co-component} of $G$ is a subgraph of $G$ induced by the vertex set of a (connected) component of $\overline{G}$.

Let  $P_\ell$  and  $C_\ell$  denote, respectively, the path and cycle of order $\ell$.
Obviously,  $P_4$  is isomorphic to its complement.  This graph will play an important role in this paper.
A graph $G$ is said to be $P_4$-free if no induced subgraph of $G$ is isomorphic to $P_4$.

The complete and edgeless graphs of order
$\ell$  will be denoted by $K_\ell$ and $S_\ell$, respectively.
Clearly, these two graphs are complementary.
A {\em clique} (respectively, a {\em stable set}) of a graph
is a  a set of pairwise adjacent (respectively, non-adjacent) vertices.
The inclusion maximal cliques and stable sets in  $V$  are called {\em maximal}, while
the cliques and stable sets of the maximal cardinality
($\omega$  and  $\alpha$, respectively) are called {\em maximum}.
These numbers  $\omega = \omega(G)$  and  $\alpha = \alpha(G)$  are referred to as
the {\em clique} and {\em stability} numbers of  $G$, respectively.

A graph   $G$  is called {\it well-covered} if every
maximal stable set of it is also maximum, that is, of size  $\alpha(G)$.
These graphs are well studied in the literature;
see, for example, \cite{BH09, BH11, CS93, Plu70, Plu93, SS92, TanTar97}.

\medskip

\begin{definition}
\label{substituting-and-lex-product}
For two vertex-disjoint graphs $G$ and $H$ and a vertex $v\in V(G)$,
{\it substituting $H$ into $G$ for $v$} means deleting $v$
and joining every vertex of $H$ to those vertices of $G$ which
have been adjacent to $v$.
The resulting graph is denoted by $G_v[H]$.

The {\it lexicographic product} of graphs $G$ and $H$ is the graph $G[H]$ with
the vertex-set $V(G)\times V(H)$, where two vertices $(u,x)$ and $(v,y)$ are adjacent if and only if
either $\{u,v\}\in E(G)$ or $u = v$ and $\{x,y\}\in E(H)$ (see, e.g.,~\cite{HIK});
in other words, this graph is obtained from $G$ by substituting $H$ for every vertex of $G$.
\end{definition}

\begin{remark}
In our paper, the families of graphs closed
with respect to lexicographic product and taking the complement
will play an important role.
The classic (and non-trivial) example is provided by the family of perfect graphs.
The Berge weak perfect graph conjecture asserted that
this family is closed under taking the complement.
Lov\'asz's \cite{Lov72a} proof of this conjecture is based on the lemma
stating that the family is closed under the substitution.
Fulkerson \cite{Ful71} was very close, but failed  to prove this lemma.
\end{remark}

\subsection{Main results}
In this paper we will study CIS circulants. Some basic facts related to the circulants and CIS graphs
will be given in the next two subsections of the Introduction.
In particular, we recall that both the CIS graphs and the circulants
are closed with respect to both operations, taking the complement
and the lexicograpic product.

It is known that $\alpha(G) \omega(G)  \leq |V(G)|$
for every circulant  $G$; see Section~\ref{sec:basic-circulants}.
We will show that a circulant  $G$ is CIS if and only if
$G$  and  $\overline G$  are both well-covered
and $\alpha(G) \omega(G)  = |V(G)|$; see Section~\ref{sec:characterization-CIS-circulants}.
The simplest subfamily of the CIS graphs is formed by the $P_4$-free graphs.
The characterization of the $P_4$-free circulants is known; see Section~\ref{sec:P4-free}.
However, it appears that there are other CIS circulants.
The minimal one is of order $36$, it was found by an exhaustive computer search;
the next two are of order $60$, see Section~\ref{sec:examples}. 

For every non-negative integer $k$ we introduce
the family of  $k$-paired circulants; see Section~\ref{sec:paired}.
We show that every $P_4$-free circulant is $k$-paired for some $k$,
yet, the value of $k$ may have to be arbitrarily large; see Section~\ref{sec:P4-free}.
%
We characterize the $2$-paired CIS circulants explicitly; see Section~\ref{sec:2-paired}.
We are not aware  of any CIS circulant that cannot be obtained from the $2$-paired CIS circulants
by the operations of taking the complement and lexicographic product.

%

\subsection{Circulants}\label{sec:basic-circulants}
For a positive integer
$n$ we denote by $[n]=\{1,2,\ldots,n\}$ the set of positive integers up to $n$.
We extend this notation to $n = 0$ by setting $[0] = \emptyset$.
Given a positive integer $n$, let us consider a set of positive integers $D\subseteq [n-1]$
such that $d\in D \Leftrightarrow n-d\in D$. The {\it circulant graph}
$C_n(D)$ is the graph with vertex-set $\mathbb{Z}_n$
in which two distinct vertices $i,j\in \mathbb{Z}_n$ are adjacent if and only if
$$i-j\pmod n\in D\,.$$ Note that a { circulant} is a Cayley graph over a cyclic group.
For example, a cycle $C_\ell$ is the circulant graph $C_\ell(\{1,\ell-1\})$.
In contrast, $P_4$  is not a circulant.

Whenever we write $i+j$ for two vertices $i,j$ of a circulant $C_n(D)$,
addition is performed modulo $n$.
For a circulant $G=C_n(D)$, we write $D(G) = D$ and refer to $D$ as a {\it distance set} of $G$.

Many graph theoretic properties of circulants can be formulated in terms of arithmetic
properties of $D$ and $n$. The following three lemmas are straightforward.

\begin{lemma}\label{lemma:circulant-complements}
The family of circulant graphs is closed under taking complements.

Furthermore, if $G = C_n(D)$ then $\overline{G} = C_n(\overline D)$, where
$\overline D = [n-1] \setminus D$.
\qed
\end{lemma}

For two positive integers $a$ and $b$, we denote by $\gcd(a,b)$ the greatest common divisor of $a$ and $b$.
Similarly, $\lcm(a,b)$ denotes the least common multiple of $a$ and $b$.
This notation naturally extends to arbitrary sequences (or sets) of integers.

\begin{sloppypar}
\begin{lemma}\label{prop:circulants-connectedness}
The number $c$ of connected components of a circulant $G = C_n(D)$ is equal to \hbox{$c = \gcd(D\cup\{n\})$}.
In particular, $G$ is connected if and only if $\gcd(D\cup\{n\}) = 1$.
Moreover, every connected component of $G$ is isomorphic to the graph $C_{n/c}(D/c)$, where $D/c = \{d/c\,:\,d\in D\}$.
\qed
\end{lemma}
\end{sloppypar}

%

\begin{lemma} \label{lemma:circulants-bipartite}
A connected circulant  $G = C_n(D)$  on at least two vertices
is bipartite if and only if  $n$  is even, while every  $d \in D$  is odd.
\qed
\end{lemma}

The next claim is a little bit more complicated, but also well known;
see. e.g., \cite{GGSTU93, GT94}.

\begin{proposition}\label{alpha-omega-le-n}
For every circulant $G$, the inequality  $\alpha(G) \omega(G) \le |V(G)|$ holds.
\end{proposition}

\begin{proof}
For completeness, we give here a short proof.
Let $G$ be a circulant of order $n$, let $C$ be a maximum clique, and let
$S$ be a maximum stable set in $G$.
Fix a vertex $j_0\in S$. For every $i\in C$, let $S_i = \{i-j_0+j\mid j\in S\}$.
Notice that $|S_i| = |S|  = \alpha(G)$ for all $i\in C$.
It is enough to show that the sets $S_i$ are pairwise disjoint, since this will imply
$$n\ge \left|\bigcup_{i\in C}S_i\right|= \sum_{i \in C}|S_i| = \sum_{i \in C}\alpha(G) = \alpha(G)\omega(G)\,.$$
Suppose for a contradiction that there exist two distinct vertices $i_1,i_2\in C$ such that
$S_{i_1}\cap S_{i_2}\neq \emptyset$. Let $k\in S_{i_1}\cap S_{i_2}$.
Then, since $k\in S_{i_1}$, there exists a $j_1\in S$ such that $k = i_1-j_0+j_1$.
Similarly, $k = i_2-j_0+j_2$ for some $j_2\in S$, since $k\in S_{i_2}$.
Therefore,
$i_1-j_0+j_1=i_2-j_0+j_2$, implying
$i_1-i_2=j_2-j_1$. This however is a contradiction, since
$i_1-i_2 \pmod n\in D(G)$ (due to $\{i_1,i_2\}\in E(G))$,
while
$j_2-j_1 \pmod n\not \in D(G)$ (due to $\{j_2,j_1\}\not\in E(G))$.
\end{proof}

The next statement was shown in \cite{GGSTU93,GT94}
(perhaps, earlier) and
recently rediscovered in \cite{Hoshino-thesis}.

\begin{proposition}
\label{prop:lex-product-circulants}
For every two circulant graphs $G$ and $H$,
their lexicographic product $G[H]$ is a circulant.
In particular, if $G = C_n(D)$ and $H= C_m(F)$, then
$G[H] = C_{nm}(T)$ where
$$T = \cup_{j = 0}^{m-1}(D+jn)\cup nF$$
where $D+jn = \{d+jn\,:\,d\in D\}$ and
$nF = \{nd\,:\,d\in F\}$.
\qed
\end{proposition}

\begin{remark}
Several classes of circulants were studied in the literature, including
well-covered~\cite{BH09, BH11, Hoshino-thesis},
perfect \cite{GGSTU93},
odd hole and odd anihole free (Berge) \cite{GGSTU93},
minimal imperfect \cite{BG93, Gri84},
partitionable  \cite{BBGMP98, BG93, Deb56, deCaen90, CGPW79, Gri84, KS06},
even hole free \cite{ABG05}, and
kernel-less oriented \cite{AFG98}.
 See \cite{GV99, Hoshino-thesis} for some additional classes of circulants.
Each of the above classes of circulants
(explored by computational experiments, in many cases) appear to be
the closures of some nice and simple subclasses under certain natural operations.
However, the questions of whether the obtained closures are complete, that is, if they contain all circulants with the considered property, are typically open,
except for
the simplest ones, such as connected, bipartite, or $P_4$-free circulants.
For example, it is still not known whether the number systems suggested
in 1956 by de Bruijn \cite{Deb56} provide all partitionable circulants (see \cite{Gri84});
although such conjecture was verified for all circulants of order
at most 36 in \cite{BBGMP98} and of at most 361 in  \cite{KS06}.
\end{remark}

\subsection{CIS graphs}
\label{sec:basic-CIS}

\begin{definition}
\label{dfn:CIS}
We say that a graph  $G$  is a {\em CIS} graph or that
it has the CIS property if
every maximal clique $C$ and every maximal stable set
$S$  in  $G$ intersect, that is, $C \cap S \neq \emptyset$.
\end{definition}


Probably Berge was the first
who paid attention to this family. 
In early 90s, Chv\'atal invited his student W. Zang
to study it, 
motivated by earlier observations of Berge and Grillet; see \cite{Zang} for more details.
The name CIS (Cliques Intersect Stable sets) was suggested in \cite{ABG06}.

\begin{lemma}\label{lemma:CIS-complements}
The family of CIS graphs is closed under taking complements.
\end{lemma}

\proof
This follows immediately from the definition, since
$K_\ell$  and  $S_\ell$  are complementary.
\qed

\smallskip

The next property just a little bit more difficult;
see, for example, \cite{ABG06}, where it
is extended from graphs to a more general case of the so-called  $d$-graphs.

\begin{lemma} 
For every two graphs $G$ and $H$ and every vertex $v\in V(G)$,
the graph $G_v[H]$ is CIS if and only if $G$ and $H$ are CIS.
\qed
\end{lemma}

{\bf Sketch of the proof}.
To show this claim, one should simply compare the maximal cliques
and stable sets of  $G$  and  $H$  to those of  $G_v[H]$.
\qed

\begin{corollary}\label{cor:lex-product-CIS}
For every $G$ and $H$, the lexicographic product $G[H]$ is CIS if and only if $G$ and $H$ are CIS.
\qed
\end{corollary}

The next statement is also very simple.

\begin{proposition}
A non-connected graph $G$ is CIS if and only if every connected component of $G$ is CIS.
A not co-connected graph $G$ is CIS if and only if all of its co-components are CIS.
\end{proposition}

%

{\bf Sketch of the proof}.
One should just compare the maximal cliques
and stable sets of the graph to those of its connected components (co-components).
\qed

For the proof of the next claim
see, e.g., \cite{Gur11}, where it
is extended from graphs to $d$-graphs.

\begin{proposition}
\label{prop:P_4-free-CIS}
Every $P_4$-free graph is CIS.
\qed
\end{proposition}


Currently, no good characterization or recognition algorithm
for the CIS graphs is known. Possible reasons for this
as well as more information about CIS graphs can be found in the Appendix.

\section{A characterization of CIS circulants}\label{sec:characterization-CIS-circulants}

\begin{theorem}\label{thm:circulant-CIS-characterization}
A circulant $G$ is a CIS graph if and only if
all maximal stable sets are of size $\alpha(G)$,
all maximal cliques are of size $\omega(G)$, and $\alpha(G)\omega(G)= |V(G)|$.

In particular, both  $G$ and  $\overline G$  are well-covered
whenever  $G$  is a CIS circulant.
\end{theorem}

\begin{proof}
Let $G$ be a CIS circulant of order $n$. Let $C$ be a maximal clique, and let $S$ be a maximal stable set in $G$.
Let $c = |C|$ and $s = |S|$.
Label the vertices of $C$ with distinct labels from the set $\{1,\ldots, c\}$, to obtain a labeled clique.
Label the vertices of $S$ with distinct labels from the set $\{1,\ldots, s\}$, to obtain a labeled stable set.
Consider the $n$ rotated copies $C_0 = C$, and $C_1,\ldots, C_{n-1}$
of the labeled clique $C$, and the $n$ rotated copies $S_0 = S$, and $S_1,\ldots, S_{n-1}$ of the labeled stable set $S$.
By the circular symmetry of $G$, every $C_i$ is a maximal clique and every $S_i$ is a maximal stable set.
We will now assign pairs of labels to vertices of $G$, as follows.
For every $i\in \mathbb{Z}_n$ and every $j\in \mathbb{Z}_n$, clique $C_i$ and stable set $S_j$ intersect
in a unique vertex $v_{ij}\in V(G)$. We assign to $v = v_{ij}$ the pair $(\ell_1,\ell_2)$ where
$\ell_1$ is the label of $v$ in $C_i$ and $\ell_2$ the label of $v$ in $S_j$.

Denote by ${\cal C}\times {\cal S}$ the set of all pairs $\{(C_i,S_j)\mid 1\le i,j\le n\}$, and
by $L = \{1,\ldots, c\}\times \{1,\ldots, s\}$ the set of all label pairs.
For a vertex $v\in V(G)$, let $L(v)$ denote the set of label pairs assigned to $v$.
Since every pair $(C_i,S_j)\in {\cal C}\times {\cal S}$ generates exactly one label pair assignment, we have
$$n^2 = |{\cal C}\times {\cal S}| = \sum_{v\in V(G)}|L(v)| \le n\cdot |L| = n|C||S|$$
and consequently
$n\le |C||S|$.
Since $|C||S|\le n$ holds for every circulant (by Proposition~\ref{alpha-omega-le-n}),
this implies that equality $|C||S| = n$ holds for every maximal clique $C$ and every maximal stable set
$S$. Choosing $C$ to be a maximum clique,
this implies that every maximal stable sets is of size $n/\omega(G)$,
and consequently every maximal stable sets is of size $\alpha(G)$, implying
$\alpha(G)\omega(G) = n$.
A symmetric argument can be used to show that every maximal clique is of size $\omega(G)$.

Conversely, suppose that $G$ is a circulant such that
all maximal stable sets are of size $\alpha(G)$,
all maximal cliques are of size $\omega(G)$, and $\alpha(G)\omega(G)= |V(G)|$.
Suppose for a contradiction that there exists a disjoint pair $(C,S)$ where $C$ is a maximal clique and
$S$ is a maximal stable set. Do the same labeling procedure as above, assigning
a label pair to a vertex in an intersection $C_i\cap S_j$ only if this intersection is nonempty.
Now, every pair $(C_i,S_j)\in {\cal C}\times {\cal S}$ generates at most one label pair assignment, and in fact
the $n$ diagonal pairs $(i,i)$ do not generate any assignment.
On the other hand, every vertex $v\in V(G)$ is assigned at least $|C||S| = \omega(G)\alpha(G)$ label pairs.
Indeed, for every $i\in C$ and every $j\in S$, the pair $(C_{i'},S_{j'})\in {\cal C}\times {\cal S}$ where
$i' = v-i\pmod n$ and $j' = v-j\pmod n$ is a pair of a clique and a stable set such that $v\in C_{i'}$ and $v\in S_{j'}$.
Hence $v$ is assigned a label pair when $C_{i'}$ and $S_{j'}$ are considered. Since
the assignments $i\mapsto i'$ and $j\mapsto j'$
are injective, we indeed have $|L(v)|\ge \alpha(G)\omega(G)$ for all $v\in V(G)$.
Putting it all together, we obtain the contradicting chain of
inequalities
$$n^2 = n\alpha(G)\omega(G) \le \sum_{v\in V(G)}|L(v)| \le |{\cal C}\times {\cal S}|-n= n^2-n\,.$$
This implies that $G$ is CIS.
\end{proof}

\begin{examples}
Let us illustrate Theorem~\ref{thm:circulant-CIS-characterization} with some examples of non-CIS
circulant graphs that violate at least one of the three conditions on the right side of the
equivalence:
\begin{enumerate}
  \item The $5$-cycle $C_5$ is a circulant graph in which all maximal stable sets are of size $\alpha(C_5) = 2$,
all maximal cliques are of size $\omega(C_5) = 2$ but $\alpha(C_5)\omega(C_5) = 4 < 5 = |V(C_5)|$.

  \item The $6$-cycle $C_6$ is a circulant graph in which
  all maximal cliques are of size $\omega(C_6) = 2$ and $\alpha(C_6)\omega(C_6) = 6  = |V(C_6)|$,
however not all maximal stable sets are of the same size.
  \item A similar example, with the role of maximal cliques and maximal sets interchanged, is given by the complement of $C_6$.
\end{enumerate}
The above examples show that none of the three conditions is implied by the other two, not even
within the class of circulants.
\end{examples}

\begin{remark}
For general (non-circulant) graphs neither of the two sides of the equivalence in Theorem~\ref{thm:circulant-CIS-characterization}
implies the other one:
\begin{itemize}
  \item The $3$-vertex path $P_3$ is a CIS graph in which not all maximal stable sets are of the same size.
  \item The $4$-vertex path $P_4$ is a graph in which all maximal stable sets are
  of the same size, all maximal cliques are of the same size, and
  $\alpha(P_4)\omega(P_4) = 4 = |V(P_4)|$. However, $P_4$ is not a CIS graph, since
  the two midpoints of it form a maximal clique $C$ that is disjoint from the maximal stable set $S$ consisting of the two endpoints of the path.
\end{itemize}
\end{remark}

%
%

\section{Paired circulants}\label{sec:paired}

\begin{sloppypar}
\begin{definition}\label{definition-a1b1a2b2}
For a non-negative integer $k$, a circulant $G = C_n(D)$ will be called {\em $k$-paired} if there
exist $k$ ordered pairs of positive integers
$(a_1, b_1),\ldots, (a_k, b_k)$ such that $a_ib_i\mid n$ for all $i\in [k]$ and
\begin{equation}\label{equation:distance-set}
D = \bigcup_{i= 1}^k D_i\,,\quad \textrm{where}\quad D_i= \Big\{d\in [n-1]\,:\,\textrm{$a_i \mid d$~~and~~$a_ib_i \centernot\mid d$}\Big\}\,.
\end{equation}
If this is the case and $k\ge 1$, we will also say that $G$ is
the {\em circulant of order $n$ generated by $a_1,b_1,\ldots, a_k,b_k$}, and denote it by
$C(n; a_1, b_1;\ldots;a_k, b_k)$.
If $k = 0$ then $D = \emptyset$, hence $G$ is edgeless; in this case we will use the notation $C(n;\emptyset)$.
A circulant that is $k$-paired for some $k$ is called {\em paired}.
\end{definition}
\end{sloppypar}


\begin{remarks}~
\begin{enumerate}
\item If $d\in D$ then $n-d\in D$, as required.
Indeed, if
$d\in [n-1]$ such that $a_i \mid d$ and $a_ib_i \centernot\mid d$
then $a_i\mid n-d$ (since
$a_i \mid n$) and
$a_ib_i \centernot\mid n- d$ (since $a_ib_i \mid n$).
\item In the definition of a $k$-paired circulant, we allow repetition of pairs and addition of the pair $(1,1)$. While neither of these operations change the graph),
    allowing them has the nice property that the classes of $k$-paired circulants form an increasing family of classes of paired circulants: if $k\le \ell$, then every $k$-paired circulant is also $\ell$-paired.



%
\end{enumerate}
\end{remarks}

\newpage
\begin{examples}~
\begin{enumerate}
  \item The edgeless graph $S_n$ of order $n$ is a
  $0$-paired circulant: $S_n = C(n;\emptyset)$.
  The complete graph $K_n$ of order $n$ is a $1$-paired circulant:  $K_n = C(n;1,n)$.
%
%

\item For every $1$-paired circulant $G = C(n;a,b)$, its complement $\overline{G}$ is $2$-paired. Indeed,
$d\in D(\overline G)$ if and only if either $a\centernot \mid d$ or $ab\mid d$. Hence,
 $\overline{G} = C(n;1,a;ab,n)$.

 Furthermore, for every $1$-paired circulant
 $G = C(n;a,b)$ with $a = 1$ or $b = 1$, its complement $\overline{G}$ is also $1$-paired.
 More specifically, if $G = C(n;1,b)$ then
 $\overline{G} = C(n;b,n/b)$, while if $G = C(n;a,1)$ then $G$ is edgeless and $\overline{G} = C(n;1,n)$ (for example).

  \item A cycle $C_n$ of order $n$ is a paired circulant if and only if $n\in \{3,4,6\}$.
  Paired circulant representations of $C_3$, $C_4$ and $C_6$ are:  $$C_3 = C(3;1,3),\quad\quad C_4 = C(4;1,2),\quad\quad C_6 = C(6; 1,2; 1,3)\,.$$
\end{enumerate}
\end{examples}

\subsection{CIS paired circulants}

The family of paired circulants is a good source of CIS circulants.  Our first infinite family of CIS paired circulants is given by the $1$-paired circulants, generalizing the complete and the edgeless graphs (which are obviously CIS).

\begin{theorem}\label{thm:1-paired-CIS}
Every $1$-paired circulant is CIS.
\end{theorem}

A proof of Theorem~\ref{thm:1-paired-CIS} will be given in Section~\ref{sec:P4-free}.

 The following theorem shows that the problem of characterizing CIS $k$-paired circulants
can be reduced to the case of connected and co-connected $k$-paired circulants.

\begin{theorem}\label{thm:CIS-k-paired-connected-co-connected}
Let $G$ be a $k$-paired circulant of order $n$ generated by $a_1,b_1,\ldots, a_k,b_k$.
Then:
\begin{enumerate}
  \item[$(i)$] If $G$ is not connected, then
  $b_i>1$ for some $i\in [k]$ and $d = \gcd(\{a_i\,:\,i\in [k]\textrm{ and } b_i>1\}) > 1$. Furthermore,
  $G$ is CIS if and only if the $k$-paired circulant of order $n/d$ generated by $a_1/d,b_1,\ldots, a_k/d,b_k$ is CIS.
  \item[$(ii)$] If $G$ is not co-connected, then there exists some $\ell\in [k]$ such $a_\ell = 1$.
  Moreover, for each such integers $\ell$, graph $G$ is CIS if and only if either $k = 1$ or
$k\ge 2$ and the
  $(k-1)$-paired circulant of order $n/b_\ell$ generated by
$A_1,B_1,\ldots, A_{\ell-1},B_{\ell-1}, A_{\ell+1}, B_{\ell+1},\ldots, A_k,B_k$ is CIS, where for all $i\in [k]\setminus\{\ell\}$ we have
$$A_i = \frac{a_i}{\gcd(a_i, b_\ell)}\quad \textrm{ and }\quad
B_i = b_i\cdot \frac{\gcd(a_i, b_\ell)}{\gcd(a_ib_i, b_\ell)}\,.$$
\end{enumerate}
\end{theorem}

A proof of Theorem~\ref{thm:CIS-k-paired-connected-co-connected} will be given at the end of Section~\ref{sec:k-paired-properties}.

For $k = 2$, we complete the characterization of CIS $k$-paired circulants in the following theorem.

\begin{theorem}\label{thm:2-paired-CIS}
Let $G$ be a connected and co-connected $2$-paired circulant of order $n$ generated by $a_1,b_1,a_2,b_2$.
Then, $G$ is CIS if and only if $\gcd(a_1b_1,a_2b_2) = 1$.
\end{theorem}

A proof of Theorem~\ref{thm:2-paired-CIS} will be given in Section~\ref{sec:2-paired}.
%
%
%
%
%

\subsection{Examples of CIS and non-CIS paired circulants}\label{sec:examples}

We now give some concrete examples of CIS and non-CIS paired circulant graphs.
In order to describe their maximal cliques and maximal stable sets, the following notation will be useful. Let $G$ be a circulant of order $n$.
Given a sequence of positive integers $\sigma = (d_1,\ldots, d_r)$,
such that
$\sum_{i = 1}^rd_i = n$,
we say that a set of vertices $X\subseteq V(G)$ {\em is generated by} $\sigma$ if there exists a vertex $i\in V(G)$ such that
$$X = \bigg\{i + \sum_{j = 1}^pd_j\,:\,1\le p\le r\bigg\}\,.$$
Notice that $i\in X$ since $\sum_{i = 1}^rd_i = n$ and additions are performed modulo $n$.
For a positive integer $p$ and an arbitrary sequence of positive integers $\sigma$, we denote by $\sigma^p$
the sequence obtained by concatenating $p$ copies of $\sigma$. More formally, if $\sigma = (d_1,\ldots, d_r)$, then
$$\sigma^p = (d_1^{(1)},\ldots, d_r^{(1)}, \ldots, d_1^{(p)},\ldots, d_r^{(p)})\,,$$
where $d^{(i)}_j = d_j$ for all $i\in [p]$.
The lists of maximal stable sets and maximal cliques for examples below
were obtained with the help of the code MACE (MAximal Clique Enumerator, ver.~2.0) for generation of all maximal cliques of a graph due to
Takeaki Uno~\cite{Uno-Mace}.

\begin{enumerate}
\label{example-2-paired-not-CIS}
   \item Let $G$ be the $2$-paired circulant $C(12;2,2;3,2)$. By Propositions~\ref{prop:k-paired-connected} and~\ref{prop:k-paired-co-connected}, respectively, $G$ is connected and co-connected. According to Theorem~\ref{thm:2-paired-CIS}, $G$ is not CIS.
  Its distance set is
    \begin{eqnarray*}
    D &=& D_1\cup D_2 \\
     &=& \{2,6,10\}\cup\{3,9\}\\
      &=& \{2,3,6,9,10\}\,.
  \end{eqnarray*}
  Every maximal clique of $G$ is generated by some sequence from the set $$\{(2, 10), (3)^4\}\,.$$
  Maximal cliques are of two different sizes, namely
$$\frac{b_1b_2\gcd(a_2,a_1b_1)}{\gcd(a_1b_1,a_2b_2)} = 2\quad \textrm{ and } \quad \frac{b_1b_2\gcd(a_1,a_2b_2)}{\gcd(a_1b_1,a_2b_2)} = 4$$
  (cf.~Proposition~\ref{prop:a1-a2-cliques} on p.~\pageref{prop:a1-a2-cliques}
  and Proposition~\ref{prop:maximal-a_1-cliques} on p.~\pageref{prop:maximal-a_1-cliques}).

  Every maximal stable set of $G$ is generated by some sequence from the set $$\{(1, 4, 7), (1, 7, 4), (4)^3\}\,.$$
   All maximal stable sets are of size~$3$.

  An example of a pair $(C,S)$ such that $C\cap S = \emptyset$
  where $C$ is a maximal clique and $S$ a maximal stable set of $G$ is given by
  $C = \{0,2\}$ and $S = \{1,5,9\}$.

  \item Let $G$ be the $2$-paired circulant $C(36;2,2;3,3)$. By Propositions~\ref{prop:k-paired-connected} and~\ref{prop:k-paired-co-connected}, respectively, $G$ is connected and co-connected. According to Theorem~\ref{thm:2-paired-CIS}, $G$  is also CIS.
  Its distance set is
  \begin{eqnarray*}
    D &=& D_1\cup D_2 \\
     &=& \{2,6,10,14,18,22,26,30,34\}\cup\{3,6,12,15,21,24,30,33\}\\
     & =& \{2,3,6,10,12,14,15,18,21,22,24,26,30,33,34\}\,.
  \end{eqnarray*}
  Every maximal clique of $G$ is generated by some sequence from the set $$\{(2, 10)^3, (3, 3, 12)^2,(6)^6\}\,.$$
  All maximal cliques are of size $b_1b_2 = 6$
  (cf.~Proposition~\ref{prop:maximal-a_1-cliques} on p.~\pageref{prop:maximal-a_1-cliques}).

  Every maximal stable set of $G$ is generated by some sequence
  from the set $$\{(1, 4, 4, 19, 4, 4), (1, 7, 1, 8, 11, 8), (4, 5, 4, 7, 9, 7)\}\,.$$
  All maximal stable sets are of size $a_1a_2 = 6$
  (cf.~Proposition~\ref{prop:well-covered} on p.~\pageref{prop:well-covered}).

  \item The $2$-paired circulant $G = C(60;2,2;3,5)$ of order $60$ is a connected and co-connected CIS circulant.
  Its distance set is
  \begin{eqnarray*}
    D &=& D_1\cup D_2 \\
     &=& \{2,6,10,14,18,22,26,30,34,38,42,46,50,54,58\}\\
     && \cup\{3,6,12,15,21,24,30,33,39,42,48,51,57\}\\
     &=& \{2,3,6,10,12,14,15,18,21,22,24,26,30,33,34,38,39,42,46,48,50,51,54,57,58\}\,.
  \end{eqnarray*}

  Every maximal clique of $G$ is generated by some sequence from the set
  $$\begin{array}{l}
  \{(2, 10)^5,(3, 3, 3, 3, 18)^2, (3, 3, 6, 12, 6)^2, (3, 6, 3, 9, 9)^2, (6)^{10}\}\,.
\end{array}$$
  All maximal cliques are of size $b_1b_2 = 10$.

  Every maximal stable set of $G$ is generated by some sequence from the set
$$\begin{array}{l}
  \{(1, 4, 11, 4, 25, 15), (1, 7, 8, 29, 8, 7), (1, 15, 1, 15, 13, 15), (1, 15, 25, 4, 11, 4),\\
  ~\,(4, 4, 7, 4, 4, 37), (4, 11, 4, 13, 15, 13),
(5, 8, 7, 8, 17, 15), (5, 15)^3, (5, 15, 17, 8, 7, 8)\}\,.
\end{array}$$
   All maximal stable sets are of size $a_1a_2 = 6$.

  \item Another example of a connected and co-connected CIS circulant on $60$ vertices is given by $2$-paired circulant $G=C(60;2,2;5,3)$.
  Its distance set is
  \begin{eqnarray*}
    D &=& D_1\cup D_2 \\
     &=& \{2,6,10,14,18,22,26,30,34,38,42,46,50,54,58\}\cup\{5,10,20,25,35,40,50,55\}\\
     &=& \{2,5,6,10,14,18,20,22,25,26,30,34,35,38,40,42,46,50,55,54,58\}\,.
  \end{eqnarray*}

  Every maximal clique of $G$ is generated by some sequence from the set
$$\{(2, 18)^3, (5, 5, 20)^2, (6, 14)^3, (10)^6\}\,.$$
  All maximal cliques are of size $b_1b_2 = 6$.

  Every maximal stable set of $G$ is generated by some sequence from the set
$$\begin{array}{l}
\{(1, 3, 4, 4, 4, 29, 4, 4, 4, 3),
(1, 3, 4, 8, 1, 15, 13, 4, 4, 7),
(1, 3, 4, 8, 21, 8, 4, 3, 1, 7),\\
~\,(1, 3, 8, 1, 3, 12, 17, 4, 8, 3),
(1, 3, 8, 4, 17, 12, 3, 1, 8, 3),
(1, 3, 9, 3, 1, 11, 4, 13, 4, 11),\\
~\,(1, 3, 9, 3, 12, 9, 8, 4, 3, 8),
(1, 3, 12, 1, 15, 1, 12, 3, 1, 11),
(1, 7, 1, 7, 1, 7, 8, 13, 8, 7),\\
~\,(1, 7, 4, 4, 13, 15, 1, 8, 4, 3),
(1, 7, 8, 1, 12, 3, 12, 1, 8, 7),
(1, 8, 3, 4, 8, 9, 12, 3, 9, 3),\\
~\,(1, 8, 4, 3, 8, 4, 9, 8, 7, 8),
(1, 8, 7, 8, 9, 4, 8, 3, 4, 8),
(3, 4, 4, 4, 9, 15, 9, 4, 4, 4),\\
~\,(4, 4, 7, 4, 4, 9, 4, 11, 4, 9),
(3, 9)^5\}\,.
\end{array}$$
   All maximal stable sets are of size $a_1a_2 = 10$.

  \item By Corollary~\ref{cor:lex-product-CIS}, CIS graphs are closed under lexicographic product, and, by Proposition~\ref{prop:lex-product-paired-circulants} (on p.~\pageref{prop:lex-product-paired-circulants}), so are the paired circulants.
      Therefore, CIS paired circulants are also closed under lexicographic product.
      For example,
      the lexicographic product $G = H[H]$, where $H$ is the CIS $2$-paired circulant $C(36;2,2;3,3)$, is
      a CIS $4$-paired circulant $C(1296; 2,2;3,3;72,2;108,3)$.
      This also shows that there exist CIS $k$-paired circulants for arbitrarily large $k$.
\end{enumerate}

\subsection{Properties of paired circulants}\label{sec:k-paired-properties}

In the rest of this section, we prove some results for general $k$-paired circulants.

\begin{proposition}\label{prop:lex-product-paired-circulants}
The family of paired circulants is closed under the lexicographic product.
More specifically, if
$G = C(n;a_1,b_1; \ldots; a_k, b_k)$
and
$H = C(m;a_1',b_1';\ldots; a_\ell', b_\ell')$,
then,
 $$G[H] \cong C(nm;a_1,b_1;\ldots; a_k, b_k; na_1',b_1'; \ldots; na_\ell', b_\ell')\,.$$
\end{proposition}

\begin{proof}
Let $G = C(n;a_1,b_1; \ldots; a_k, b_k)$
and
$H = C(m;a_1',b_1';\ldots; a_\ell', b_\ell')$.
Denoting by $D$ and $F$ the distance sets of $G$ and $H$, respectively, we have
$$D =
\bigcup_{i = 1}^k\Big\{d\in [n-1]\,:\,a_i\mid d \textrm{ and } a_ib_i \centernot\mid d\Big\}$$
and
$$F =
\bigcup_{j = 1}^\ell\Big\{f\in [m-1]\,:\,a_j'\mid f\textrm{ and } a_j'b_j' \centernot\mid f\Big\}\,.$$

\begin{sloppypar}
By Proposition~\ref{prop:lex-product-circulants}, we have
$G[H] = C_{nm}(T)$ where
$T = \cup_{j = 0}^{m-1}(D+jn)\cup nF\,.$
To establish the proposition, we will show that
$T$ is equal to the distance set $T'$ of  the paired circulant
$C(nm;a_1,b_1;\ldots; a_k, b_k; na_1',b_1'; \ldots; na_\ell', b_\ell')$, which is given by the expression
$$T' =
\bigcup_{i = 1}^k\Big\{t\in [nm-1]\,:\,a_i\mid t \textrm{ and } a_ib_i \centernot\mid t\Big\}
\cup
\bigcup_{j = 1}^\ell\Big\{t\in [nm-1]\,:\,na_j' \mid t \textrm{ and } na_j'b_j' \centernot\mid t\Big\}\,.$$
\end{sloppypar}

First, let $t\in T$. Then, either there exists an integer $j\in [0,m-1]$ such that $t \in D+jn$ or $t\in nF$.

In the former case, $t = d+jn$ for some $d\in D$ and $j\in [0,m-1]$. Let $i\in [k]$ be an integer such that
$a_i\mid d$ and $a_ib_i \centernot\mid d$.
Since $a_ib_i \mid n$, we have that $a_i\mid t$ and $a_ib_i \centernot\mid t$.
Since $j\in [0,m-1]$ and $d\in [n-1]$, it follows that $t = d+jn\in [mn-1]$. Hence $t\in T'$.

In the latter case, $t = nf$ for some $f\in F$. Let $j\in [\ell]$ be an integer such that
$a_j'\mid f$ and $a_j'b_j' \centernot\mid f$.
Since $a_j'\mid f$ and $a_j'b_j' \mid m$, we have that $na_j'\mid t$ and $na_j'b_j' \centernot\mid t$.
Since $f\in [m-1]$, it follows that $t\in [mn-1]$. Hence $t\in T'$.

Second, let $t\in T'$. Then, $t\in[nm-1]$, and
either there exists an integer $i\in [k]$ such that
$a_i\mid t$ and $a_ib_i \centernot\mid t$, or
there exists an integer $j\in [\ell]$ such that
$na_j'\mid t$ and $na_j'b_j' \centernot\mid t$.

In the former case, let $i\in [k]$ be an integer such that
$a_i\mid t$ and $a_ib_i \centernot\mid t$.
Let $d = t \pmod{n}$ and $j = (t-d)/n$.
If $d = 0$ then $n$ divides $t$, which is impossible since $a_ib_i \centernot\mid t$.
Hence, $d\in [n-1]$ and $a_i\mid d$ and $a_ib_i \centernot\mid d$ (since $a_ib_i \mid n$).
Since $d\le t$, integer $j$ is non-negative. Moreover, $j\le m-1$,
since otherwise we would obtain a contradicting chain of inequalities
$nm\le t-d\le t\le nm-1$. This implies that $t = d+jn$ for some $d\in D$ and $j\in [0,m-1]$, hence $t\in T$.

In the latter case, let $j\in [\ell]$ be an integer such that
$na_j'\mid t$ and $na_j'b_j' \centernot\mid t$.
Let $f = t/n$. Since $na_j'\mid t$, it follows that $f$ is an integer.
Clearly, $f\ge 0$, and also $f\le m-1$, since otherwise $t\ge nm$.
Hence, $f\in [0,m-1]$. Furthermore, the definition of $f$ together with the properties
$na_j'\mid t$ and $na_j'b_j' \centernot\mid t$ imply that
$a_j'\mid f$ and $a_j'b_j' \centernot\mid f$.
Consequently, $f\in F$, and $t \in nF\subseteq T$.
\end{proof}

\begin{proposition}\label{prop:reduction-lex-product}
Suppose that $G$ is the circulant of order $n$ generated by $a_1,b_1, \ldots, a_k, b_k$ where $k\ge 1$.
Let $d = \lcm(a_1b_1,\ldots, a_kb_k)$.
Then, $G$ is isomorphic to the lexicographic product of
$C(d; a_1,b_1;\ldots; a_k,b_k)$ and $S_{\frac{n}{d}}$.
\end{proposition}

\begin{proof}
First, recall that the edgeless graph
$S_{\frac{n}{d}}$ is isomorphic to the $0$-paired circulant $C(\frac{n}{d}; \emptyset)$.
By Proposition~\ref{prop:lex-product-paired-circulants}, the lexicographic product of
paired circulants
$C(d; a_1,b_1;\ldots; a_k,b_k)$ and $S_{\frac{n}{d}} = C(\frac{n}{d}; \emptyset)$
is isomorphic to the $k$-paired circulant
$C(n; a_1,b_1;\ldots; a_k,b_k) =G$.
\end{proof}

\begin{sloppypar}
\begin{proposition}\label{prop:k-paired-connected}
Let $G$ be a $k$-paired circulant $G = C(n;a_1,b_1;\ldots; a_k,b_k)$.
Then, the number $d$ of connected components of $G$
is equal to $$d=\left\{
                                     \begin{array}{ll}
                                       n, & \hbox{if $b_i = 1$ for all $i\in [k]$;} \\
                                       \gcd(A), & \hbox{otherwise;} \\
                                     \end{array}
                                   \right.$$
                                   where $A = \{a_i\,:\,i\in [k]\textrm{ and } b_i>1\}$.
In particular, $G$ is connected if and only if
either $n = 1$ or
$b_i>1$ for some $i\in [k]$ and $\gcd(A) = 1$.
Furthermore, $G$ is isomorphic to the lexicographic product of the edgeless graph
$S_d$ and the $k$-paired circulant
$C(\frac{n}{d}; \frac{a_1}{d},b_1;\ldots;\frac{a_1}{d},b_k)$.
\end{proposition}
\end{sloppypar}

%

\begin{proof}
First, let us show that the number of connected components of $G$ is indeed given by the above expression.
By Lemma~\ref{prop:circulants-connectedness},
it is enough to show that $d = \gcd(D\cup\{n\})$, where $D$ is the distance set of $G$ given by~\eqref{equation:distance-set}.
If $b_i = 1$ for all $i\in [k]$ (in particular, this is trivially the case if $k = 0$),
then $G$ is edgeless and $d = n$, as specified by the expression.
Suppose now that $b_i>1$ for some $i\in [k]$.
Since $b_i>1$ for every $a_i\in A$, we have $a_i\in D$. Therefore, $A\subseteq D\cup\{n\}$, and every common divisor of $D\cup\{n\}$
is also a common divisor of $A$, which shows that $\gcd(A)\ge \gcd(D\cup\{n\})$.
On the other hand, the definition of $k$-paired circulants implies that
every common divisor of $A$ is also a common divisor of $D\cup\{n\}$, which
shows that $\gcd(D\cup\{n\})\ge \gcd(A)$.

To prove the last part of the proposition, we proceed
as in the proof of Proposition~\ref{prop:reduction-lex-product}.
The edgeless graph
$S_{d}$ is isomorphic to the $0$-paired circulant $C(d;\emptyset)$. Moreover, by
Proposition~\ref{prop:lex-product-paired-circulants}, the lexicographic product of
paired circulants
$S_{d} = C(d;\emptyset)$ and $C(\frac{n}{d}; \frac{a_1}{d},b_1;\ldots;\frac{a_1}{d},b_k)$
is isomorphic to the $k$-paired circulant
$C(n; a_1,b_1;\ldots; a_k,b_k) =G$.
\end{proof}

To prove Proposition~\ref{prop:k-paired-co-connected} below, we will need the following straightforward observation
relating the operations of the lexicographic product and the complement.

\begin{proposition}\label{prop:complement-and-lex-product}
For every two graphs $G$ and $H$, the graphs $\overline{G[H]}$ and $\overline{G}[\overline H]$ are isomorphic.\qed
\end{proposition}

\begin{sloppypar}
\begin{proposition}\label{prop:k-paired-co-connected}
A $k$-paired circulant $G = C(n;a_1,b_1;\ldots; a_k,b_k)$
is co-connected if and only if
either $n = 1$ or
$a_i\ge 2$ for all $i\in [k]$. Furthermore, if there exists an $\ell\in [k]$ such that
$a_\ell = 1$, then for every such $\ell$,
  graph $G$ is isomorphic to the lexicographic product
of the complete graph $K_{b_\ell}$ and the $(k-1)$-paired circulant of order $n/b_\ell$ generated
  by $A_1,B_1,\ldots, A_{\ell-1},B_{\ell-1}, A_{\ell+1}, B_{\ell+1},\ldots,A_k,B_k$
   such that for all $i\in [k]\setminus\{\ell\}$ we have
$$A_i = \frac{a_i}{\gcd(a_i, b_\ell)}\quad \textrm{ and }\quad
B_i = b_i\cdot \frac{\gcd(a_i, b_\ell)}{\gcd(a_ib_i, b_\ell)}
\,.$$
\end{proposition}
\end{sloppypar}

\begin{proof}
Let us first show the first part of the proposition, that is,
that $G$ is co-connected if and only if
either $n = 1$ or
$a_i\ge 2$ for all $i\in [k]$.

The case $n = 1$ is trivial, so let $n\ge 2$.

Suppose first that $G$ is co-connected, and suppose
for a contradiction that $a_\ell = 1$ for some $\ell\in [k]$.
By definition of $D$, if $b_\ell\centernot\mid d$ and $d\in [n-1]$, then
$d\in D$.
Hence, every distance in the complementary distance set $\overline D = [n-1]\setminus D$
is divisible by $b_\ell$. In particular, $\gcd(\overline D)\ge b_\ell>1$, and Lemma~\ref{prop:circulants-connectedness}
and Lemma~\ref{lemma:circulant-complements} imply that
the complementary circulant $\overline G$ is not connected, contrary to the assumption that $G$ is co-connected.

Suppose now that $a_i\ge 2$ for all $i\in [k]$. Then, by the definition of $D$ we have
$1\in \overline D$, which, by Lemmas~\ref{lemma:circulant-complements} and~\ref{prop:circulants-connectedness}
implies that the complementary circulant $\overline G$ is connected.

To prove the last part of the proposition, suppose that $\ell\in [k]$ is such that $a_\ell = 1$.
For simplicity, let us assume that $\ell = 1$.

First, we handle the case when $k = 1$.
In this case $G = C(n;1,b_1)$, and its complement $\overline{G}$ is
the $1$-paired circulant $G = C(n; b_1, n/b_1)$.
Hence, by Proposition~\ref{prop:k-paired-connected}
the graph $\overline G$ is isomorphic
to the lexicographic product of the edgeless graph
$S_{b_1}$ and
the $1$-paired circulant
$C(n/b_1; 1,n/b_1)$.
Consequently, since
$C(n/b_1; 1,n/b_1)$ is the complete graph of order $n/b_1$, Proposition~\ref{prop:complement-and-lex-product} implies that
graph $G$ is isomorphic to the lexicographic product of the complete graph
$K_{b_1}$ and  the (edgeless) $0$-paired circulant
of order $n/b_1$. This establishes the proof for the case $k = 1$.

Now, suppose that $k\ge 2$.
We will show that the distance set of $\overline{G}$ is equal to the distance set of the
lexicographic product of graphs $S_{b_1}$ and $\overline{G'}$, where $G'$ is the $(k-1)$-paired circulant defined in the proposition.
Since both $\overline{G}$  and $S_{b_1}[\overline{G'}]$ are circulant graphs on $n$ vertices, the claim
will then follow from Proposition~\ref{prop:complement-and-lex-product}.
We have the following:
\begin{itemize}
  \item The distance set of the graph $G$
is equal to
$$D = \Big\{d\in [n-1]\,:\,\textrm{$b_1 \centernot\mid d$}\Big\}\cup \bigcup_{i= 2}^k \Big\{d\in [n-1]\,:\,\textrm{$a_i \mid d$~~and~~$a_ib_i \centernot\mid d$}\Big\}\,.$$

  \item The distance set of the graph $\overline G$
is equal to
\begin{eqnarray*}
T &=& [n-1]\setminus D\\
 &=& \Big\{t\in [n-1]\,:\,\textrm{$b_1 \mid t$}\Big\}
\cap \bigcap_{i= 2}^k \Big\{t\in [n-1]\,:\,\textrm{$a_i \centernot\mid t$~~or~~$a_ib_i \mid t$}\Big\}\\
 &=& \Big\{b_1t'\,:\,t'\in [n/b_1-1]\textrm{~~and~~($\forall i\in \{2,\ldots, k\}$)($a_i \centernot\mid b_1t'$~~or~~$a_ib_i \mid b_1t'$)}\Big\}\,.
\end{eqnarray*}

  \item The distance set of the graph
$G' = C(n/b_1; A_2,B_2;\ldots;A_k,B_k)$
is equal to
$$D' = \bigcup_{i= 2}^k \Big\{d'\in [n/b_1-1]\,:\,\textrm{$A_i \mid d'$~~and~~$A_iB_i \centernot\mid d'$}\Big\}\,.$$
  \item The distance set of the graph
$\overline{G'}$
is equal to
\begin{eqnarray*}
T' &=& [n/b_1-1]\setminus D'\\
&=& \bigcap_{i= 2}^k \Big\{t'\in [n/b_1-1]\,:\,\textrm{$A_i \centernot\mid t'$~~or~~$A_iB_i \mid t'$}\Big\}\,.
\end{eqnarray*}

\item
By Proposition~\ref{prop:lex-product-circulants}, the distance set of the lexicographic product of
graphs $S_{b_1}$ and $\overline{G'}$ is equal to
\begin{eqnarray*}
b_1\cdot T' &=& b_1\cdot \bigg(\bigcap_{i= 2}^k \Big\{t'\in [n/b_1-1]\,:\,\textrm{$A_i \centernot\mid t'$~~or~~$A_iB_i \mid t'$}\Big\}\bigg)\\
 &=& \Big\{b_1t'\,:\,t'\in [n/b_1-1]\textrm{~~and~~($\forall i\in \{2,\ldots, k\}$)($A_i \centernot\mid t'$~~or~~$A_iB_i \mid t'$)}\Big\}\\
 &=& \Big\{b_1t'\,:\,t'\in [n/b_1-1]\textrm{~~and~~($\forall i\in \{2,\ldots, k\}$)($a_i \centernot\mid b_1t'$~~or~~$a_ib_i \mid b_1t'$)}\Big\}\\
 &=& T\,.
\end{eqnarray*}
Only the third equality above requires some justification.
The equality follows from the following two equivalences:
\begin{eqnarray}
\label{eq:Ai1}A_i \mid t' &\textrm{~~if and only if~~}&a_i \mid b_1t'\\
\label{eq:Ai2}
A_iB_i \mid t' &\textrm{~~if and only if~~}& a_ib_i \mid b_1t'
\end{eqnarray}
Let us verify these two equivalences.  For~\eqref{eq:Ai1}, observe that,
on the one hand, if $A_i\mid t'$ then there exists an integer $r$
such that $$t' = r A_i = \frac{r a_i}{\gcd(a_i, b_1)}\,,$$ therefore
$$b_1 t' = \frac{b_1 r a_i}{\gcd(a_i, b_1)} = a_i r'\,,$$
where $$r' = \frac{b_1 r}{\gcd(a_i, b_1)}$$ is integer,
and consequently $a_i \mid b_1t'$.
On the other  hand, if  $a_i \mid b_1t'$, then there exists an integer $r$
such that
$$b_1t' =  r a_i  = r\, {\gcd(a_i, b_1)} A_i\,.$$
Since $A_i$ and $b_1$ are relatively prime, this implies that
$A_i\mid t'$.

Equivalence~\eqref{eq:Ai2} can be proved similarly, using the fact that
$$A_iB_i = \frac{a_ib_i}{\gcd(a_ib_i, b_1)}\,.$$
 \end{itemize}
\end{proof}

\medskip

\noindent{\bf Proof of Theorem~\ref{thm:CIS-k-paired-connected-co-connected}}.

\begin{proof}[Proof (Theorem~\ref{thm:CIS-k-paired-connected-co-connected})]
Part $(i)$ of the theorem follows from Proposition~\ref{prop:k-paired-connected} and
Corollary~\ref{cor:lex-product-CIS}. Similarly, part $(ii)$ follows from
Lemma~\ref{lemma:CIS-complements}, Proposition~\ref{prop:k-paired-co-connected} and Corollary~\ref{cor:lex-product-CIS}.
\end{proof}

\section{CIS $2$-paired circulants}\label{sec:2-paired}

In this section, we prove Theorem~\ref{thm:2-paired-CIS}.
The theorem will be derived in Section~\ref{sec:proof} from the results of the previous sections and of the rest of this section.
More specifically, in Sections~\ref{sec:cliques} and~\ref{sec:stable-sets} we will analyze the structure of maximal cliques and maximal stable sets in
$2$-paired circulants, respectively.

\subsection{Maximal cliques}\label{sec:cliques}

Let $G$ be a $2$-paired circulant of order $n$ generated by $a_1, b_1, a_2, b_2$.
We say that a clique $C$ in $G$ is an {\it $a_1$-clique} if
$i\equiv j \pmod {a_1}$ holds for every two vertices $i,j\in C$.
Similarly, a clique $C$ is said to be an {\it $a_2$-clique} if
$i\equiv j \pmod {a_2}$ holds for every two vertices $i,j\in C$.

\begin{proposition}\label{prop:a1-a2-cliques}
Every clique $C$ in the graph $C(n;a_1,b_1;a_2,b_2)$ is either an $a_1$-clique or an $a_2$-clique.
\end{proposition}

\begin{proof}
Let $C$ be a clique in $G$. It follows directly from Definition~\ref{definition-a1b1a2b2} that
every two vertices $i,j\in C$ satisfy either
$i\equiv j \pmod {a_1}$ or
$i\equiv j \pmod {a_2}$ (or both).
Suppose that $C$ is not an $a_1$-clique.
Then, the relation of congruence modulo $a_1$ has at least two equivalence classes
$C_1,\ldots, C_r$. We claim that in this case, every two vertices $i,j\in C$ are congruent modulo $a_2$.
Indeed, if $i\not\equiv j\pmod{a_1}$ then, as observed above, this implies $i\equiv j \pmod {a_2}$.
On the other hand, if $i\equiv j\pmod{a_1}$, then $i$ and $j$ belong to the same equivalence class $C_{p}$.
Let $k$ be an arbitrary vertex from an equivalence class $C_{p'}$ such that $p'\neq p$.
Then, $i\equiv k\pmod{a_2}$ and
\hbox{$k\equiv j\pmod{a_2}$}, and thus,
since the relation of congruence modulo $a_2$ is transitive,
we infer that
\hbox{$i\equiv j\pmod{a_2}$} holds as well.
Thus, $C$ is an $a_2$-clique in this case.
\end{proof}

\begin{sloppypar}
\begin{proposition}\label{prop:maximal-a_1-cliques}
Every maximal $a_1$-clique in the graph
$C(n;a_1,b_1;a_2,b_2)$ such that $n = \lcm(a_1b_1,a_2b_2)$
is of size exactly $$\frac{b_1b_2\gcd(a_2,a_1b_1)}{\gcd(a_1b_1,a_2b_2)}\,.$$
Every maximal $a_2$-clique in the graph
$C(n;a_1,b_1;a_2,b_2)$
such that $n = \lcm(a_1b_1,a_2b_2)$
is of size exactly
$$\frac{b_1b_2\gcd(a_1,a_2b_2)}{\gcd(a_1b_1,a_2b_2)}\,.$$
\end{proposition}
\end{sloppypar}

\begin{proof}
Let $C$ be an $a_1$-clique in the graph $C(n;a_1,b_1;a_2,b_2)$.
Due to the circular symmetry of $G$, we may assume that $0\in C$.
Hence, every vertex $i\in C$ can be written in a unique way as
$$i = a_1\Bigg(r_i+\bigg(\alpha_i+t_i\frac{a_2}{\gcd(a_2,a_1b_1)}\bigg)b_1\Bigg)$$
for some $r_i\in [0,b_1-1]$, $\alpha_i\in [0,\frac{a_2}{\gcd(a_2,a_1b_1)}-1]$,
and $t_i\in [0, \frac{b_2\gcd(a_2,a_1b_1)}{\gcd(a_1b_1,a_2b_2)}-1]$.

We claim that if $r_i = r_j$ for some $i,j\in C$, then
$\alpha_i =  \alpha_j$.
Indeed, suppose that $r_i = r_j$ but
(say) $\alpha_i > \alpha_j$.
Then
\begin{eqnarray*}
i-j
&=&
a_1b_1\bigg((\alpha_i-\alpha_j)+(t_i-t_j)\frac{a_2}{\gcd(a_2,a_1b_1)}\bigg)\\
&=&
a_1b_1\bigg((\alpha_i-\alpha_j)+(t_i-t_j)\frac{a_2\lcm(a_2,a_1b_1)}{a_1b_1a_2}\bigg)\\
&=&
a_1b_1(\alpha_i-\alpha_j)+(t_i-t_j)\lcm(a_2,a_1b_1)\,.
\end{eqnarray*}
Hence, $i\equiv j\pmod{a_1b_1}$.
Moreover, since
$a_1b_1(\alpha_i-\alpha_j)<a_1b_1\frac{a_2}{\gcd(a_2,a_1b_1)} = \lcm(a_2,a_1b_1)$, the
definition of
the least common multiple implies
that $a_2 \centernot\mid a_1b_1(\alpha_i-\alpha_j)$. Consequently, $i\not \equiv j\pmod{a_2}$, which
contradicts the fact that $i$ and $j$ are adjacent.

The above observation implies that for every $i\in C$, the value of $\alpha_i$ is uniquely determined with the value of
$r_i$. Thus, $\alpha_i$ is a function of $r_i$, and we write $\alpha_i = \alpha(r_i)$.
Consequently, if $r_i = r_j$ and $t_i= t_j$ for some $i,j\in C$, then $i = j$.

Therefore, for every $r\in[0,b_1-1]$ there exists
at most one $\alpha_r\in [0, \frac{a_2}{\gcd(a_2,a_1b_1}-1]$
such that there exists a vertex $i\in C$ with
$r_i = r$ and $\alpha_i = \alpha_r$.
Moreover, for every such pair $(r,\alpha_r)$ and every
$t\in [0,
\frac{b_2\gcd(a_2,a_1b_1)}{\gcd(a_1b_1,a_2b_2)}-1]$,
there is at most one vertex $i\in C$ such that
$r_i = r$, $\alpha_i = \alpha_r$, and
$t_i = t$.
Hence, the total number of vertices in $C$ is at most
$$b_1\cdot \frac{b_2\gcd(a_2,a_1b_1)}{\gcd(a_1b_1,a_2b_2)} = \frac{b_1b_2\gcd(a_2,a_1b_1)}{\gcd(a_1b_1,a_2b_2)}\,.$$

To conclude the proof, suppose for a contradiction that $C$ has strictly less than
$\frac{b_1b_2\gcd(a_2,a_1b_1)}{\gcd(a_1b_1,a_2b_2)}$ vertices.
We analyze two cases.

{\em Case 1. There exists an integer $\tilde r\in[0,b_1-1]$
such that there is no vertex $i\in C$ with $r_i = \tilde r$.}

Let $\tilde i = a_1b_1\tilde r\,.$
Clearly, $\tilde i$ is a vertex of $G$, and the assumption on $\tilde r$ implies that $\tilde i\not\in C$.
We claim that $\tilde i$ is adjacent to every $i\in C$.
Indeed, for every $i\in C$ we have
\begin{eqnarray*}
\tilde i-i &=&
a_1\Bigg((\tilde r-r_i)-\bigg(\alpha_i+t_i\frac{a_2}{\gcd(a_2,a_1b_1)}\bigg)b_1\Bigg)\,.
\end{eqnarray*}
Hence,
$\tilde i\equiv i\pmod{a_1}$ but $\tilde i\not \equiv i\pmod{a_1b_1}$ and consequently $\tilde i$ is adjacent to $i$.
Since the choice of $i\in C$ was arbitrary, this contradicts the maximality of $C$.

\medskip
{\em Case 2. For every integer $r\in[0,b_1-1]$ there exists a vertex $i\in C$ with $r_i = r$.}

In this case, the above derivation of the inequality $|C|\le \frac{b_1b_2\gcd(a_2,a_1b_1)}{\gcd(a_1b_1,a_2b_2)}$ together with the assumption that the inequality is strict imply that there exist integers
$\tilde r$, $\tilde t$ with $\tilde r\in[0,b_1-1]$,
and
$\tilde t\in [0,\frac{b_2\gcd(a_2,a_1b_1)}{\gcd(a_1b_1,a_2b_2)}-1]$ such that no vertex $i\in C$ satisfies
$r_i = \tilde r$ and $t_i = \tilde t$.

Let
$$\tilde i =
a_1\Bigg(\tilde r+\bigg(\alpha(\tilde r)+\tilde t \cdot\frac{a_2}{\gcd(a_2,a_1b_1)}\bigg)b_1\Bigg)\,.$$
It is easy to verify that $\tilde i\in [n-1]$, that is, $\tilde i$ is a vertex of $G$. Moreover, by the choice of $\tilde r$ and $\tilde t$, we have
$\tilde i\not\in C$. We will reach a contradiction with maximality of $C$ by showing that vertex $\tilde i$ is adjacent to every vertex $i\in C$.

For $i\in C$ with $r_i\neq \tilde r$, we derive (similarly as in Case~$1$ above) $\tilde i\equiv i\pmod{a_1}$ and
$\tilde i\not \equiv i\pmod{a_1b_1}$; consequently $\tilde i$ is adjacent to $i$.

Suppose now that vertex $i\in C$ is such that $r_i= \tilde r$. Then $\alpha_i = \alpha(\tilde r)$ and $t_i \neq \tilde t$.
Therefore
\begin{eqnarray*}
i-\tilde i &=&  \frac{a_1b_1a_2(t_i-\tilde t)}{\gcd(a_2,a_1b_1)}\\
&=&  \lcm(a_2,a_1b_1)(t_i-\tilde t)\,.
\end{eqnarray*}
Hence,
$a_2\mid i-\tilde i$, and also
$a_1b_1\mid i-\tilde i$.
If also $a_2b_2\mid i-\tilde i$, then $\lcm(a_1b_1,a_2b_2) = n \mid i-\tilde i$, which is impossible
since $1\le |i-\tilde i|\le n-1$.
Therefore $i\equiv \tilde i\pmod{a_2}$ but
$i\not\equiv \tilde i\pmod{a_2b_2}$,
which implies that $\tilde i$ and  $i$ are adjacent.

This completes the proof of Case 2 and with it the proof of the first part of the proposition.

The second part follows by symmetry.
\end{proof}

\begin{corollary}\label{cor:co-well-covered}
All maximal cliques in the graph $C(n;a_1,b_1;a_2,b_2)$ such that $n = \lcm(a_1b_1,a_2b_2)$
are of the same size if and only~if
$\gcd(a_2,a_1b_1) = \gcd(a_1,a_2b_2)\,.$
\end{corollary}

\subsection{Maximal stable sets}\label{sec:stable-sets}

\begin{sloppypar}
Let us now consider maximal stable sets in a $2$-paired circulant $G$ of order $n$ generated by $a_1, b_1, a_2, b_2$.
such that $n = \lcm(a_1b_1,a_2b_1)$.
To every pair of distinct non-adjacent vertices $i$ and $j$ in $G$, let us associate a two dimensional label $\ell(i,j)\in\ZZ_+^2$, defined by 
$$\ell(i,j)_1 = i-j\mod{a_1}\quad\textrm{~~and~~} \quad\ell(i,j)_2 = i-j\mod{a_2}\,.$$
Pairs $(i,j)$ of distinct non-adjacent vertices of $G$ will also be referred to as {\em directed non-edges} (of~$G$).
\end{sloppypar}

\begin{proposition}\label{prop:label-(0,0)-not-in-range}
For every  directed non-edge $(i,j)$ of $G = C(n;a_1,b_1;a_2,b_2)$ such that $n = \lcm(a_1b_1,a_2b_1)$, we have
$\ell(i,j)\neq (0,0)$.
\end{proposition}

\begin{proof}
Suppose that $\ell(i,j)= (0,0)$ for a directed non-edge $(i,j)$ of $G$.
Then $i\equiv j \pmod{a_1}$ and $i\equiv j \pmod{a_2}$.
Since $i$ and $j$ are non-adjacent,
$i\equiv j \pmod{a_1}$
implies
$i\equiv j \pmod{a_1b_1}$ and similarly,
$i\equiv j \pmod{a_2}$
implies
$i\equiv j \pmod{a_2b_2}$.
Consequently, $i \equiv j\pmod{n}$, a contradiction.
\end{proof}

\begin{proposition}\label{prop:triangles}
If $(i,j)$, $(j,k)$ and $(k,i)$ are directed non-edges of $G = C(n;a_1,b_1;a_2,b_2)$
such that $n = \lcm(a_1b_1,a_2b_1)$, then
$$\ell(i,j)_1 +  \ell(j,k)_1 + \ell(k,i)_1 \equiv 0 \pmod{a_1}$$
and
$$\ell(i,j)_2 +  \ell(j,k)_2 + \ell(k,i)_2 \equiv 0 \pmod{a_2}\,.$$
\end{proposition}

\begin{proof}
By symmetry, it suffices to prove the first congruence.
We have
$$\ell(i,j)_1 +  \ell(j,k)_1 + \ell(k,i)_1 \equiv
(i-j) + (j-k) + (k-i)\equiv 0 \pmod{a_1}\,.$$
\end{proof}

\begin{proposition}\label{prop:labels-in-max-stable-set}
Let $S$ be a maximal stable set
in $G = C(n;a_1,b_1;a_2,b_2)$
such that $n = \lcm(a_1b_1,a_2b_1)$.
Suppose that $0\in S$.
Then, for all $j,j'\in S\setminus \{0\}$ such that
$j\neq j'$, we have $\ell(0,j)\neq \ell(0,j')$.
\end{proposition}

\begin{proof}
Suppose that $\ell(0,j)\neq \ell(0,j')$ for some
$j,j'\in S\setminus \{0\}$. Then
$0-j\equiv 0-j'\pmod{a_1}$ as well as
$0-j\equiv 0-j'\pmod{a_2}$. Hence  $a_1\mid j'-j$ and
$a_2\mid j'-j$.
Since $j'$ and $j$ are non-adjacent, we have $a_1b_1\mid j'-j$
and similarly $a_2b_2\mid j'-j$.
Consequently, $n\mid j'-j$, which implies $j' = j$.
\end{proof}

Propositions~\ref{prop:label-(0,0)-not-in-range} and~\ref{prop:labels-in-max-stable-set} together with the circular symmetry of $G$
imply the following.

\begin{corollary}\label{cor:upper-bound-on-alpha}
Every maximal stable set $S$ in  $G = C(n;a_1,b_1;a_2,b_2)$
such that $n = \lcm(a_1b_1,a_2b_1)$
satisfies $|S|\le a_1a_2$.
\end{corollary}

\begin{proposition}\label{prop:well-covered}
If $\gcd(a_1b_1,a_2b_2) = 1$ then every maximal stable set $S$
in $G = C(n;a_1,b_1;a_2,b_2)$
such that $n = \lcm(a_1b_1,a_2b_1)$ satisfies
$|S|= a_1a_2$.
\end{proposition}

\begin{proof}
Suppose that $\gcd(a_1b_1,a_2b_2) = 1$ and that $S$ is a maximal stable set in $G$
such that $|S|<a_1a_2$.
Due to the circular symmetry of $G$, we may assume that $0\in S$.
Let $F = \{\ell(0,j)\,:\,j\in S\setminus\{0\}\}$.
Let $(u,v)$ be an arbitrary element of the (nonempty) set
$$\Big([0, a_1-1]\times [0,a_2-1]\Big)\setminus \Big(F\cup\{(0,0)\}\Big)\,.$$

We will show that there exists a vertex $x\in V(G)\setminus S$ such that $S\cup \{x\}$ is a stable set, where
$x$ is of the form
\begin{equation}\label{eq1}
x = u + \alpha a_1 + \gamma a_1b_1 = v + \beta a_2 + \delta a_2b_2
\end{equation}
for some $\alpha, \beta, \gamma, \delta$ such that
$\alpha \in [0, b_1-1]$, $\beta \in [0, b_2-1]$,
$\gamma \in [0, a_2b_2-1]$, $\delta \in [0, a_1b_1-1]$, and, in addition,
the following conditions are met:
\begin{equation}\label{eq2}
  \textrm{if } u = 0 \textrm{ then } \alpha = 0\,,
\end{equation}
\begin{equation}\label{eq3}
  \textrm{if } v = 0 \textrm{ then } \beta = 0\,.
\end{equation}
To this end, let us consider first the following congruence:
\begin{equation}\label{eq4}
u + \alpha a_1 + \gamma a_1b_1\equiv v + \beta a_2 + \delta a_2b_2\pmod{n}\,.
\end{equation}

\begin{claim}\label{claim:gamma-delta}
For every two integers $\alpha$ and $\beta$, there exist integers
$\gamma \in [0, a_2b_2-1]$ and $\delta \in [0, a_1b_1-1]$ such that equation~(\ref{eq4}) holds.
\end{claim}

\begin{proof}
Since
$a_1b_1$ and $a_2b_2$ are relatively prime, the Diophantine
equation
$$a_1b_1\gamma'  - a_2b_2\delta' = v-u + \beta a_2-\alpha a_1$$
has a solution $(\gamma',\delta')$. Taking modulo $n$ both sides and shifting
$\gamma'$ and $\delta'$ by an appropriate multiples of
$a_2b_2$ and $a_1b_1$, respectively,
we can find $\gamma$ and $\delta$ satisfying the conditions of the claim.
\end{proof}

Let us partition the set $S\setminus\{0\}$ into three pairwise disjoint subsets $S_1$, $S_2$, $S_3$, where
$$S_1 = \{j\in S\setminus\{0\}\,:\, (\exists v')(\ell(0,j) = (u,v'))\}\,,$$
$$S_2 = \{j\in S\setminus\{0\}\,:\, (\exists u')(\ell(0,j) = (u',v))\}\,,$$
$$S_3 = \{j\in S\setminus\{0\}\,:\, \ell(0,j)_1\neq u\,, \ell(0,j)_2 \neq v\}\,.$$

\begin{claim}\label{claim:alpha}
Suppose that $S_1\neq\emptyset$. Then, there exists an integer $\alpha\in [0,b_1-1]$
and integers $\lambda_j$ for $j\in S_1$
such that for all $j\in S_1$, it holds that
$$j = u+\alpha a_1+\lambda_j a_1b_1\,.$$
Furthermore, if $u = 0$ then $\alpha = 0$.
\end{claim}

\begin{proof}
By the definition of $S_1$, every $j\in S_1$ satisfies
$j\equiv u \pmod{a_1}$, and if $u = 0$ then $j\equiv u\pmod{a_1b_1}$.
This implies the claimed form of $j$ with $\alpha$ possibly depending on $j$.
It also implies that if $u = 0$ then $\alpha = 0$.

If $|S_1| > 1$, then for every two distinct elements $j, j'\in S$
we have
$$j = u+\alpha a_1+\lambda_j a_1b_1$$
and
$$j' = u+\alpha' a_1+\lambda_{j'}a_1b_1\,,$$ which implies
$$j-j' = (\alpha-\alpha')a_1+(\lambda_j-\lambda_{j'})a_1b_1\,.$$
Since $\{j,j'\}$ is a non-edge in $G$ and since
$a_1\mid j-j'$, we must have $\alpha\equiv \alpha' \pmod{b_1}$.
\end{proof}

By symmetry, we can also show the following.

\begin{claim}\label{claim:beta}
Suppose that $S_2\neq\emptyset$. Then, there exists an integer $\beta \in [0,b_2-1]$
and integers $\mu_j$ for $j\in S_2$
such that for all $j\in S_2$, it holds that
$$j = v+\beta a_2+\mu_j a_2b_2\,.$$
Furthermore, if $v = 0$ then $\beta  = 0$.
\end{claim}

Now we are ready to define $x$.
If $S_1 \neq \emptyset$ then we set $\alpha$ according to Claim~\ref{claim:alpha}.
If $S_1 = \emptyset$ then we set $\alpha = 0$.
Analogously, if
$S_2 \neq \emptyset$ then we set $\beta$ according to Claim~\ref{claim:beta}.
If $S_2 = \emptyset$ then we set $\beta = 0$.
Finally, we set $\gamma$ and $\delta$ according to Claim~\ref{claim:gamma-delta}, and set $x$ as in
equation~\eqref{eq1}.

\begin{claim}\label{claim:0}
Vertex $x$ is not adjacent to vertex $0$.
\end{claim}

\begin{proof}
By the definition of $x$, we have
$x\equiv u\pmod{a_1}$ and
$x\equiv v\pmod{a_2}$.
If $u\neq 0$ and $v\neq 0$, then the claim is implied.

If $u= 0$ then we also have that
$x\equiv 0\pmod{a_1b_1}$ by Claim~\ref{claim:alpha} if $S_1\neq\emptyset$
and by the definition of $\alpha$ if $S_1= \emptyset$.
In this case, $v \neq 0$ (since $(u,v)\neq (0,0)$), therefore
$x\not \equiv 0\pmod{a_2}$, proving the claim.

Analogously, Claim~\ref{claim:beta} and the definition of $\beta$ imply that
$x$ is non-adjacent to $0$ if $v = 0$.
\end{proof}

\begin{claim}\label{claim:S_1}
For every $j\in S_1$, vertex $x$ is not adjacent to vertex $j$.
\end{claim}

\begin{proof}
Let $j\in S_1$.

By the choice of $(u,v)$, we have $\ell(0,j)_2 = v'\neq v$.
Since $x\equiv v\pmod{a_2}$, we have that $a_2 \centernot\mid x-j$.

Let us also note that by Claim~\ref{claim:alpha}
and by the definition of $x$, we have
$x-j = (\gamma -\lambda_j)a_1b_1$, hence $a_1b_1\mid x-j$, proving the claim.
\end{proof}

An analogous proof shows the following.

\begin{claim}\label{claim:S_2}
For every $j\in S_2$, vertex $x$ is not adjacent to vertex $j$.
\end{claim}

\begin{claim}\label{claim:S_3}
For every $j\in S_3$, vertex $x$ is not adjacent to vertex $j$.
\end{claim}

\begin{proof}
Let $j\in S_3$. By the definition of $S_3$, we have
$\ell(0,j)_1 = u'\neq u$ and
$\ell(0,j)_2 = v'\neq v$, therefore
$a_1 \centernot\mid x-j$ and $a_2 \centernot\mid x-j$.
\end{proof}

Since by the choice of $(u,v)$, vertex $x$ cannot belong to $S$, the above
claims imply that $S$ cannot be a maximal stable set, which proves the statement
of the proposition.
\end{proof}

\begin{proposition}\label{prop:small-stable-sets}
Suppose that $a_1>1$, $a_2>1$, $\gcd(a_1,a_2b_2) = \gcd(a_2,a_1b_1) = 1$
but \hbox{$\gcd(b_1,b_2) > 1$}. Then, the graph
$G = C(n;a_1,b_1;a_2,b_2)$
such that $n = \lcm(a_1b_1,a_2b_1)$
has a stable set $S'$ of size~$3$ such that
for all stable sets $S$ with $S'\subseteq S$, it holds that
$|S|<a_1a_2$.
\end{proposition}

\begin{proof}
Without loss of generality, we may assume that $a_2\ge 3$.
Let $\alpha = a_2\cdot \gcd(b_1,b_2)$, and let $\beta$ be an integer such that
\begin{equation}\label{eq5}
2+\beta a_2\equiv 0\pmod{a_1b_1}\,.
\end{equation}
Such a solution exists in the range $[0,a_1b_1-1]$ since
$\gcd(a_2,a_1b_1)  = 1$.

Let $i = 1+\alpha a_1$ and $j = 2+\beta a_2$.
By the definition, we have $i\neq 0$ and $j\neq 0$, moreover
$i\neq j$ since $i\equiv 1\pmod{a_2}$ (as $a_2\mid \alpha$) while
$j\equiv 2\pmod{a_2}$. We claim that the set $S' = \{0,i,j\}$ is a stable set:
\begin{itemize}
  \item We have
  $i\not\equiv 0 \pmod {a_1}$, and the definition of $\alpha$ implies that
$i\not\equiv 0 \pmod {a_2}$. Therefore vertex  $i$ is not adjacent to vertex $0$.

  \item We have $j\equiv 0 \pmod {a_1b_1}$ and
  $j\equiv 2\not\equiv 0\pmod{a_2}$. Therefore vertex $j$ is not adjacent to vertex $0$.

  \item   We have $i\equiv 1\pmod{a_1}$ and
$j\equiv 0\pmod{a_1}$, hence $i\not \equiv j  \pmod{a_1}$. Moreover,
$i\equiv 1\pmod{a_2}$ and
$j\equiv 2\pmod{a_2}$, hence $i\not \equiv j  \pmod{a_2}$.
Therefore, vertices $i$ and $j$ are non-adjacent.
\end{itemize}

Let $S$ be a maximal stable set in $G$ such that $S'\subseteq S$.
We will show that $|S|<a_1a_2$, which will establish the statement of the proposition.
Suppose for a contradiction that $|S| = a_1a_2$.
By Proposition~\ref{prop:labels-in-max-stable-set},
there exists a vertex $x\in S$
such that $\ell(0,x) = (1,2)$.
By the definition of $\ell(0,x)$, we have
$x\equiv 1\pmod{a_1}$ and $x\equiv 2\pmod{a_2}$.
These congruences and the definitions of $i$ and $j$ imply that
$x\equiv i\pmod{a_1}$ and $x\equiv j\pmod{a_2}$.
Since $x$ is non-adjacent to both $i$ and $j$, we must have
$x\equiv i\pmod{a_1b_1}$ and $x\equiv j\pmod{a_2b_2}$.
This implies the existence of integers
$\lambda$ and $\mu$ such that
$x = i+\lambda a_1b_1$
and
$x = j+\mu a_2b_2$.
Therefore,
\begin{equation}\label{eq6}
j-i = \lambda a_1b_1-\mu a_2b_2\,.
\end{equation}
We have
 $i\equiv 1\pmod{\gcd(b_1,b_2)}$ by the definition of $i$, while
$j\equiv 0\pmod{\gcd(b_1,b_2)}$ by the definition of $j$ and congruence~\eqref{eq5}.
Therefore, since $\gcd(b_1,b_2)>1$, we have
$j-i \not\equiv 0 \gcd(b_1,b_2)$. On the other hand,
$\gcd(b_1,b_2)$ divides both terms in the right hand side of
equation~\eqref{eq6}, hence $\lambda a_1b_1-\mu a_2b_2\equiv 0\pmod{\gcd(b_1,b_2)}$.

This contradiction shows that our assumption about the
size of $S$ was incorrect.
\end{proof}

\subsection{Proof of Theorem~\ref{thm:2-paired-CIS}}\label{sec:proof}

Recall that Theorem~\ref{thm:2-paired-CIS} states that
if $G$ is a $2$-paired connected and co-connected circulant of order $n$ generated by $a_1,b_1,a_2,b_2$,
then $G$ is CIS if and only if $\gcd(a_1b_1,a_2b_2) = 1$.

\begin{proof}[Proof (Theorem~\ref{thm:2-paired-CIS})]~

If $n = 1$ then $a_1 = b_1 = a_2 = b_2 = 1$ and the statement of the theorem clearly holds.

So let $n\ge 2$. By Proposition~\ref{prop:k-paired-connected}, the fact that
$G$ is connected implies that $b_1\ge 2$ or $b_2\ge 2$.
If one of $b_1$ and $b_2$ is equal to $1$, say $b_2 = 1$, then $b_1>1$ and
Proposition~\ref{prop:k-paired-connected} implies that $a_1 = 1$, contrary to the fact that $G$ is co-connected and
Proposition~\ref{prop:k-paired-co-connected}. Thereofore, $b_1\ge 2$ and  $b_2\ge 2$, and, since
$G$ is connected, $\gcd(a_1,a_2) = 1$ by Proposition~\ref{prop:k-paired-connected}.

Since $G$ is co-connected, $a_1\ge 2$ and $a_2\ge 2$ by Proposition~\ref{prop:k-paired-co-connected}.

Now, let us argue that it suffices to prove the theorem for the case when
$n = \lcm(a_1b_1,a_2b_2)$. Indeed, if $d = \lcm(a_1b_1,a_2b_2)<n$, then
Proposition~\ref{prop:reduction-lex-product} implies that
$G$ is isomorphic to the lexicographic product of
the $2$-paired circulant $G' = C(d; a_1,b_1; a_2,b_2)$
and the edgeless graph $S_{\frac{n}{d}}$.
Assume that
$G'$ is CIS if and only if $\gcd(a_1b_1,a_2b_2) = 1$.
Then,  Corollary~\ref{cor:lex-product-CIS} and the fact that
$S_{\frac{n}{d}}$ is CIS imply that $G$ is CIS if and only if $G'$ is CIS.
Therefore,
$G$ is CIS if and only if $\gcd(a_1b_1,a_2b_2) = 1$.

Let us now assume that $n = \lcm(a_1b_1,a_2b_2)$. We will now verify both implications of the equivalence.

For the forward direction, assume that $G$ be CIS. Assume indirectly that
$\gcd(a_1b_1,a_2b_2) > 1$. Since $G$ is CIS, Theorem~\ref{thm:circulant-CIS-characterization}
implies that all maximal cliques of $G$ are of the same size.
Therefore, since
$n = \lcm(a_1b_1,a_2b_2)$, the condition in Corollary~\ref{cor:co-well-covered} holds, and since $\gcd(a_1,a_2) = 1$, the condition
can be simplified to $\gcd(a_2,b_1) = \gcd(a_1,b_2)$.

Since $\gcd(a_1,a_2) = 1$, it must be the case that either $\gcd(a_1,b_2) = \gcd(a_2,b_1) > 1$ or $\gcd(b_1,b_2) > 1$.
Suppose first that $\gcd(a_1,b_2) = \gcd(a_2,b_1) > 1$.
Let $d = \gcd(a_1,b_2)$. Then, $d$ divides $a_1$ as well as $a_2$, which
contradicts $gcd(a_1,a_2)=1$.

Suppose now that $\gcd(a_1,b_2) =\gcd(a_2,b_1) = 1$. Then
$\gcd(a_1b_1,a_2b_2) = \gcd(b_1,b_2) >1$.
On the one hand, Proposition~\ref{prop:maximal-a_1-cliques} implies that
all maximal cliques in $G$ are of size
$$\omega(G) = \frac{b_1b_2\gcd(a_2,b_1)}{\gcd(a_1b_1,a_2b_2)} = \frac{b_1b_2}{\gcd(b_1,b_2)} = \lcm(b_1,b_2)\,.$$
On the other hand, Proposition~\ref{prop:small-stable-sets}
implies that $G$ has a stable set $S'$ of size~$3$ such that
for all stable sets $S$ with $S'\subseteq S$, it holds that $|S|<a_1a_2$.
Since $G$ is CIS,  Theorem~\ref{thm:circulant-CIS-characterization} implies that all maximal stable sets of $G$ are of the same size
and consequently $\alpha(G) <a_1a_2$.
Thus,
$$\alpha(G)\omega(G) < a_1a_2\lcm(b_1,b_2) = \lcm(a_1b_1,a_2b_2) = n\,,$$
contrary to Theorem~\ref{thm:circulant-CIS-characterization}.

For the converse direction, suppose that $\gcd(a_1b_1,a_2b_2) = 1$.
By Proposition~\ref{prop:maximal-a_1-cliques}, all maximal cliques are of size
$\omega(G) = b_1b_2$. By Proposition~\ref{prop:well-covered}, all maximal stable sets are of size
$\alpha(G) = a_1a_2$. Consequently, $$\alpha(G)\omega(G) = a_1b_1a_2b_2 = \gcd(a_1b_1,a_2b_2)\lcm(a_1b_1,a_2b_2) = n\,,$$
and Theorem~\ref{thm:circulant-CIS-characterization} implies that $G$ is CIS.

This completes the proof of Theorem~\ref{thm:2-paired-CIS}.
\end{proof}

\section{$P_4$-free circulants}\label{sec:P4-free}

Let us recall the following well-known
characterization of the $P_4$-free graphs
and its corollaries.

\begin{proposition}[\cite{CLB81, Gur77, Gur84, Sum71}]\label{thm:P_4-free-decomposition}
A graph $G$ is $P_4$-free if and only if for every induced subgraph
$F$ of $G$ with at least two vertices, either $F$ or its complement is not connected.
\end{proposition}

\begin{proposition}[see, e.g., \cite{GGSTU93}]\label{prop:P_4-free-lex-product}
The class of $P_4$-free graphs is closed under lexicographic product.
\end{proposition}

\begin{proposition}\label{thm:1-paired-P4-free}
Every $1$-paired circulant is $P_4$-free.
\end{proposition}

\begin{proof}
We will show the theorem by induction on the number of vertices.
For $n = 1$, the statement is trivially true.

Let $G=C(n;a,b)$ be a $1$-paired circulant on $n>1$ vertices, and suppose that
the statement of the theorem holds for all graphs on less than $n$ vertices.
If $b = 1$ then $G$ is edgeless, and hence $P_4$-free.
So let $b\ge 2$.

If $a = 1$ then $G$ is not co-connected, and
by Proposition~\ref{prop:k-paired-co-connected} its complement
$\overline{G}$ is isomorphic to the lexicographic product
of the edgeless graph $S_b$ of order $b$
and the complement of the $0$-paired circulant of order $n/b$, that is,
the complete graph $K_{n/b}$ of order $n/b$.
Since edgeless and complete graphs are $P_4$-free, so is
$\overline{G}$, by Proposition~\ref{prop:P_4-free-lex-product}.
Since the graph $P_4$ is isomorphic to its complement, the
$P_4$-free graphs are also closed under taking the complement
and hence $G$ is $P_4$-free as well.


If $a>1$ then $G$ is not connected, and by Proposition~\ref{prop:k-paired-connected},
$G$ is isomorphic to the lexicographic product of the edgeless graph
$S_a$ and the $1$-paired circulant $C(\frac{n}{a}; 1,b)$.
By induction, the $1$-paied circulant $C(\frac{n}{a}; 1,b)$ is $P_4$-free.
Hence, $G$ is $P_4$-free by Proposition~\ref{prop:P_4-free-lex-product}.
\end{proof}

Propositions~\ref{thm:1-paired-P4-free} and \ref{prop:P_4-free-CIS} imply Theorem~\ref{thm:1-paired-CIS}.

\begin{theorem-1-paired-CIS}
Every $1$-paired circulant is CIS.\qed
\end{theorem-1-paired-CIS}

\begin{remark}
Theorem~\ref{thm:1-paired-CIS} and Proposition \ref{thm:1-paired-P4-free} are best possible, in the sense that not every {\it $2$-paired} circulant is CIS.
For example, the $2$-paired circulant $C(12;2,2;3,2)$ is connected (by Proposition~\ref{prop:k-paired-connected}) and also co-connected (by Proposition~\ref{prop:k-paired-co-connected}). However, by Theorem~\ref{thm:2-paired-CIS}, it is not CIS (cf.~Example~$1$ on p.~\pageref{example-2-paired-not-CIS}).
\end{remark}

\begin{theorem}\label{thm:P4-free-paired}
Every $P_4$-free circulant is paired.
\end{theorem}

\begin{proof}
We will show the theorem by induction on the number of vertices.
For $n = 1$, the statement of the theorem clearly holds.

Let $G$ be a $P_4$-free circulant on $n>1$ vertices, and suppose that
the statement of the theorem holds for all graphs on less than $n$ vertices.
By Proposition~\ref{thm:P_4-free-decomposition},
either $G$ or its complement is not connected.
Suppose first that $G$ is not connected.
Then, by Lemma~\ref{prop:circulants-connectedness}, $G$ has exactly \hbox{$d = \gcd(D\cup\{n\})$} connected components,
where $D$ is a distance set of $G$, and every connected component of $G$ is isomorphic to
$C_{n/d}(D/d)$. The graph $C_{n/d}(D/d)$ is a $P_4$-circulant, hence, by induction,
$C_{n/d}(D/d)$ is $k$-paired for some $k$, say
$C_{n/d}(D/d) = C(n/d;a_1,b_1;\ldots; a_k,b_k)$.
Since $G$ is the lexicographic product of $S_d = C(d;\emptyset)$ and
$C(n/d;a_1,b_1;\ldots; a_k,b_k)$, Proposition~\ref{prop:lex-product-paired-circulants}
implies that $G$ is isomorphic to the $k$-paired circulant $C(n;da_1,b_1;\ldots, da_k,b_k)$.

Suppose now that the complement of $G$ is not connected. Then, again
by Lemma~\ref{prop:circulants-connectedness}, there exists an integer $d>1$ such that $\overline{G}$ has exactly $d$ connected components, each of which is isomorphic to some circulant $H$. Since the complement of $H$ is a $P_4$-free circulant on less than $n$ vertices, the inductive hypothesis implies that
$\overline H$ is $k$-paired for some $k$, that is, that $\overline H = C(n/d;a_1,b_1;\ldots, a_k,b_k)$ for some positive integers $a_1,b_1,\ldots, a_k,b_k$.
Since $G$ is isomorphic to the lexicographic product of the complete graph $K_d = C(d;1,d)$ and $\overline H = C(n/d;a_1,b_1;\ldots, a_k,b_k)$.
Proposition~\ref{prop:lex-product-paired-circulants}
implies that $G$ is isomorphic to the $(k+1)$-paired circulant $C(n;1,d;da_1,b_1;\ldots, da_k,b_k)$.
\end{proof}

\begin{remark}
Theorem~\ref{thm:P4-free-paired} shows that every $P_4$-free circulant is paired.
The converse is not true, as shown by the $2$-paired circulant $C(36;2,2;3,3)$.
This motivates the following question:
Given a paired circulant $G$, how can we determine whether $G$ is $P_4$-free?
Propositions~\ref{prop:k-paired-connected} and~\ref{prop:k-paired-co-connected}
provide a recursive decomposition procedure of a given paired circulant $G$ into connected components of
$G$ or its complement. Proposition~\ref{thm:P_4-free-decomposition} implies that this procedure gives
an efficient way of determining whether a given paired circulant is $P_4$-free:
A paired circulant $G$ is $P_4$-free if and only if $G$ can be decomposed into copies of
$1$-vertex paired circulant $C(1;\emptyset)$.
\end{remark}

\begin{proposition}\label{thm:P4-free-k-paired}
For every $k$, there exists a $P_4$-free circulant that is not $k$-paired.
\end{proposition}

\begin{sloppypar}
\begin{proof}
Let $p_1, p_2, p_3, \ldots$ be an enumeration of all primes.
For every positive integer $n$, let $Q_n =K_{p_{2n-1}}[S_{p_{2n}}]$, that is,
$Q_n$ is the lexicographic product of the complete graph of order $p_{2n-1}$ and the edgeless graph of order $p_{2n}$.
Notice that $Q_n$ is of order $q_n=p_{2n-1}p_{2n}$. Let us define a sequence of circulants $\{G_n\}_{n\ge 1}$ recursively as follows:
\begin{itemize}
  \item $G_1 = Q_1$, and
  \item for $n\ge 2$, let $G_n = Q_n[G_{n-1}]$ be the lexicographic product of $Q_n$ with $G_{n-1}$.
\end{itemize}
An induction on $n$ together with Proposition~\ref{prop:P_4-free-lex-product} implies that
every $G_n$ is $P_4$-free.
Induction on $n$ and Proposition~\ref{prop:lex-product-circulants} show that
$G_n = C_{g_n}\big(D^{(n)}\big)$ where $g_n = \prod_{i = 1}^nq_i$
and the distance set $D^{(n)}$ can be computed recursively using the formulas
$$D^{(1)} = [g_1]\setminus p_{1}[p_{2}]$$
and
\begin{equation}\label{eq:disj-union}
D^{(n)}
= \Big([g_n]\setminus p_{2n-1}[g_{n-1}p_{2n}]\Big)\cup q_nD^{(n-1)}\,.
\end{equation}
It follows from the above formulas that $G_n$ is an $n$-paired circulant of order $g_n$
generated by $(a_1,b_1),\ldots, (a_n,b_n)$ where for each $i\in [n]$,
we have $a_i = \prod_{j = i+1}^nq_j$ (with $a_n = 1$) and
$b_i = p_{2i-1}$.

We will prove by induction on $n$ that for every $n\ge 1$, graph $G_n$ is not $(n-1)$-paired.

For $n =1$, the fact that $G_1$ is not $0$-paired follows from the fact that the only $0$-paired circulants are the edgeless ones, and $G_1$ is not edgeless.
Now, let $n\ge 2$, and suppose inductively that graph $G_{n-1}$ is not $(n-2)$-paired.

To show that $G_n$ is not $(n-1)$-paired, it is sufficient to show that
the distance set $D^{(n)}$ cannot be represented as the union
\begin{equation}\label{eq:union}
D^{(n)} = \bigcup_{i= 1}^{p} D_i\,,\quad \textrm{where}\quad D_i= \Big\{d\in [g_n-1]\,:\,\textrm{$\alpha_i \mid d$~~and~~$\alpha_i\beta_i \centernot\mid d$}\Big\}
\end{equation}
for some positive integers $\alpha_1,\beta_1,\ldots, \alpha_{p},\beta_{p}$ where $p\le n-1$.
Indeed, the result of the proposition will then follow by applying a result of Muzychuk~\cite{Muzychuk}
stating that if $N$ is a positive integer not divisible by the square of any prime number,
then any two circulants of order $N$ that are isomorphic have the property that their distance sets $D$ and $D'$
satisfy $D' = qD$ where $q\in \{1,2,\ldots, N-1\}$ such that $\gcd(q,N) = 1$.
In our case, we have $N = g_n$, which by construction is not divisible by the square of any prime.
Moreover, if $D' = qD$ for some $q$ as above, then it is easy to verify that if
$C_N(D)$ is a $k$-paired circulant generated by $(a_1,b_1),\ldots, (a_k,b_k)$ then so is $C_N(D')$. (In fact, $D' = D$.)

Suppose for a contradiction that $D^{(n)}$ can be represented as the union as in~\eqref{eq:union} with $p\le n-1$.
Among all such representations, take one with minimum $p$.
Since $1\in D^{(n)}$ and $p_{2n-1}\not\in D^{(n)}$, there exists an $i\in [p]$ such that $\alpha_i = 1$ and $\beta_i= p_{2n-1}$.
Without loss of generality we may assume that $\alpha_p = 1$ and $\beta_p=p_{2n-1}$.
Consequently, $D_p$ contains all distances in $D^{(n)}$ that are not divisible by $p_{2n-1}$.
By the minimality of $p$, all other $D_i$'s contain distances divisible by
$p_{2n-1}$. In fact, since all distances in $D^{(n)}$ that are
divisible by $p_{2n-1}$
are also divisible by $q_{n} = p_{2n-1}p_{2n}$, every $\alpha_i$ for $i\in \{1,\ldots, p-1\}$
can be written in the form $\alpha_i = q_n\alpha_i'$ for some positive integer $\alpha_i'$.
By equation~\eqref{eq:disj-union},
$D^{(n)}$ is the disjoint union of $D_p$ and $q_nD^{(n-1)}$.
This implies that
$$D^{(n-1)} = \bigcup_{i= 1}^{p-1} D_i'\,,\quad \textrm{where}\quad D_i'= \Big\{d\in [g_{n-1}-1]\,:\,\textrm{$\alpha_i' \mid d$~~and~~$\alpha_i'\beta_i \centernot\mid d$}\Big\}\,,$$
that is, that the graph $G_{n-1}$ is a $(p-1)$-paired circulant (generated by $(\alpha_1',\beta_1),\ldots, (\alpha_{p-1}',\beta_{p-1})$). This is a contradiction with the fact that $p-1\le n-2$ and the inductive hypothesis that $G_{n-1}$ is not $(n-2)$-paired.
\end{proof}
\end{sloppypar}

\section{Open questions and problems}

%

It is not known whether the following statements are true or false:
\begin{itemize}
 \item
 Every CIS circulant can be obtained from the $2$-paired CIS circulants by taking the complements and lexicographic products.

 \item For every CIS circulant  $G$,  either $G$  or its complement  $\overline G$  is paired.

 \item Isomorphic circulants  $C_n(D)$  and  $C_n(D')$
 either both are  $k$-paired or both are not, for any fixed  $k$.
     Clearly, this conjecture holds for the so-called {\em Caley isomorphisms},
     $D \rightarrow i D \pmod n$, where $\gcd(i, n) = 1$.
    However, for some isomorphic pairs there exist other isomorphisms \cite{Muzychuk}.

\item For (i) paired, (ii) $k$-pared, and (iii) CIS circulants
there exist only Caley isomorphisms.

\end{itemize}

\noindent
The following questions are also open:

 \begin{itemize}
 \item
 Which $k$-paired circulants are CIS?
 The answer is known only for  $k\le 2$.

 \item
 How difficult is it to determine whether a given circulant $C_n(D)$ is
 (i) paired?  (ii)  $k$-paired?  (iii) CIS?

Let us remark that the recognition problem of well-covered circulants is co-NP-complete~\cite{BH11}.

%

\end{itemize}

Another research direction is extending the results obtained in this paper
to Cayley graphs of other groups.

\subsection*{Acknowledgements}

The smallest nontrivial connected and co-connected CIS circulant is
$C(36;2,2;3,3)$, using notation of Definition~\ref{definition-a1b1a2b2}.
We found this circulant with an exhaustive search through all small circulants, with the help of
the code MACE (MAximal Clique Enumerator, ver.~2.0) for generation of all maximal cliques of a graph due to
Takeaki Uno~\cite{Uno-Mace}.
We are grateful to Jernej Vi\v ci\v c for
letting us use his computer and for help with technical matters related with running the programs,
and Tine \v Sukljan for helpful suggestions regarding implementation.

\section*{Appendix: More about  CIS graphs}

Let us notice that the graph  $P_4$  itself is not a CIS graph. Indeed, let
$E(P_4) = \{v_1 v_1', v_1' v_2', v_2' v_2\}$, then
$S = \{v_1, v_2\}$  and  $C = \{v'_1, v'_2\}$  are
disjoint maximal stable set and maximal clique, $C \cap S = \emptyset$.

However,  $P_4$  is an induced subgraph of a CIS graph
$A$ of order $5$  defined by the edge-set
$E(A) = \{v_1 v_1', v_1' v_2', v_2' v_2, v_0 v'_1, v_0 v'_2\}$.
It is easily seen that  $A$  is a CIS graph and, by construction,  $P_4$
is the subgraph of  $A$  induced by
$V(P_4) = \{v_1, v_1', v_2', v_2\} = V(A) \setminus \{v_0\}$.
Graph $A$  is called the {\em bull-graph} or the $A$-{\em graph}.

This simple example shows that
the family of CIS graphs is not hereditary, that is,
not closed under taking induced subgraphs.
Moreover, the following observation was proved in \cite{ABG06}.

\begin{proposition}
\label{proposition:each-graph-is-a-subgraph-of-a-CIS}
Every graph  $G$  is an induced subgraph of a CIS graph  $G'$.
\end{proposition}

\proof
Let us extend every maximal clique  $C$  of $G$  by a new
({\em simplicial}) vertex  $v_C$  that is connected in  $G'$
to all vertices of  $C$  and to no other.
By construction, the obtained graph  $G'$  contains  $G$
as an induced subgraph. Furthermore, it is easy to verify that  $G'$ is CIS.
\qed

\medskip

Let us make the following remarks:

\smallskip
\begin{itemize}
\item We  do not need to extend  $C$  whenever it already has a simplicial vertex in  $G$.
In particular, such an ``economical'' extension of  $P_4$
results exactly in the $A$-graph.
\item The order of  $G'$  may be exponential in the order of  $G$.
\item One can get another CIS extension of  $G$ by complementing an extension of its complement.
\end{itemize}

\medskip

The above CIS extension of $P_4$  can be generalized as follows.
A $k$-{\em comb} is a graph  $B_k$  with  $2k$  vertices,
$V(B_k) = \{v_1, \ldots, v_k; v'_1, \ldots, v'_k\}$, and
the edge-set  $E(B_k)$  such that
$S = \{v_1, \ldots, v_k\}$  is a stable set,
$C = \{v'_1, \ldots, v'_k\}$  is a  clique,
$\{v_i, v'_i  \mid  i \in [k] = \{1, \ldots, k\}\}$  is a matching,
and there are no moore edges, that is, $|E(B_k)| = \binom {k} {2} + k = k(k+1)/2$.
Let us extend  $B_k$  adding to it one vertex  $v_0$  and  $k$  edges
$v_0 v'_i$  for all $i \in [k]$, and denote the obtained graph by  $D_k$.
It is easily seen that  $D_k$  is CIS and, by construction, $B_k$
is the subgraph of  $D_k$  induced by  $V(D_k) \setminus \{v_0\} = V(B_k)$.
Graph  $D_k$  is called a {\em settled} $k$-{\em comb}.

The complementary graphs  $\overline B_k$  and $\overline D_k$  are called
an {\em anticomb} and {\em settled anticomb}, respectively.
For example, $P_4$  is a $2$-comb and $k$-anticomb simultaneously.
However, for  $k > 2$  the $k$-comb and $k$-anticomb are not isomorphic.

Let us also notice that in a  $k$-comb, as well as in a $k$-anticomb,
all  $\ell$-combs and  $\ell$-anticombs are settled, for every
$\ell < k$.

The following condition, obviously, is necessary for the CIS property to hold;
see, e.g., \cite{ABG06}.

\begin{proposition}
Every induced $k$-comb (respectively, $k$-anticomb) of a CIS graph  $G$
is an induced subgraph of a settled induced $(k+1)$-comb
(respectively, $(k+1)$-anticomb) of  $G$;
in other words, every  comb and anticomb must be settled in $G$.
\qed
\end{proposition}

For   $k=2$   this observation means that
in a CIS graph each induced $P_4$ must be settled by an  $A$-graph.
Probably, Berge was the first who noticed it in 70s; see \cite{Zang} for  more details.

\begin{example}
In 1994, Holzman demonstrated  that
the above condition is only necessary but not sufficient
for the CIS property to hold.
Let us consider $\binom {5} {2} + \binom {5} {1} = 10 + 5 = 15$  vertices
$$V = \{v_{ij}, v_k \mid i,j,k \in [5] = \{1,2,3,4,5\}, i \neq j\},$$
and define the edge-set  $E$  such that the first ten vertices,
$\{v_{ij} \mid i,j \in [5], i \neq j\}$, form a clique  $C$, the last five
$\{v_k \mid k \in [5]\}$  form a stable set  $S$,  and
a pair   $v_{ij} v_k$  is an edge if and only if   $k \in \{i,j\}$.

It is not difficult to verify that
in the obtained graph  $H$  all combs  and anticombs are settled.
For example, the $3$-anticomb induced by
$\{v_{12}, v_{13}, v_{23}, v_1, v_2, v_3\}$  is settled by  $v_{45}$, while
the $4$-comb induced by
$\{v_{12}, v_{13}, v_{14}, v_{14}, v_2, v_3, v_4, v_5\}$  is settled by  $v_1$.
Clearly,  $H$  contains no  $5$-combs or  $4$-anticombs.
It is also easy to check that every  $2$-comb in $H$   is settled.
However, $H$  is not CIS, since
the clique $C$  and stable set  $S$  are maximal and $C \cap S = \emptyset$;
see \cite{ABG06} for  more details.
\end{example}

The following theorem provides sufficient conditions for the CIS property to hold.
It was conjectured by Chvatal in 90s and
proved in \cite{DengLZ04, DengLZ05} and then
independently in \cite{ABG06}. Both proofs are lengthy and technical.

\begin{theorem}
A graph $G$ is CIS if it contains no induced
$3$-combs and anticombs and every induced $2$-comb in it is settled.
\qed
\end{theorem}

It is not known whether the following weaker conditions
are still sufficient for the CIS property to hold for a graph $G$:
(i) all induced $2$-combs, $3$-combs, and $3$-anticombs are settled and
there are no induced $4$-combs and $4$-anticombs in $G$;
(ii) all induced combs and anticombs are settled
and there is no induced Holzman graph  $H$ in $G$.

\smallskip

Thus, currently, no good characterization or recognition algorithm
for the CIS graphs is known.
One can notice certain similarity to the perfect graphs, replacing
the combs and anticombs by the odd holes and antiholes.
Yet, unlike the CIS graphs, perfect graphs form a hereditary class.

\smallskip

In contrast, the next class admits a very simple
characterization (which is not easy to prove, yet).
A graph is called {\em almost CIS} if every its
maximal clique $C$  and maximal stable set $S$  intersect, except for
a unique pair, $C_0$  and  $S_0$.
Somewhat surprisingly, the next characterization hods.

\begin{theorem}
Graph  $G = (V,E)$  is almost CIS if and only if
$V = C_0 \cup S_0$, where
$C_0$ is a maximal clique,
$S_0$ is a maximal stable set,
and  $C_0 \cap S_0 = \emptyset$;
or, in other words,  if and only if  $G$  is a  split graph with
a unique split partition. 
\end{theorem}

This claim was conjectured and
some partial results obtained in \cite{BGZ09}.
Then it was proved in \cite{WZZ}.
In particular, this theorem implies that every split graph
is either CIS or almost CIS.


\begin{thebibliography}{99}

\bibitem{ABG05}
D. Andrade, E. Boros, and V. Gurvich,
Even hole free and balanced circulants,
RUTCOR Research Report  RRR-08-2005, Rutgers University.

\bibitem{ABG06}
D.V. Andrade, E. Boros, and V. Gurvich, On Graphs Whose Maximal Cliques and Stable
Sets Intersect, RRR 17-2006, RUTCOR Research Reports, Rutgers University.

\bibitem{AFG98}
A. Apartsin, E. Ferapontova, and V. Gurvich.
A circular graph -- counterexample to the Duchet kernel conjecture,
Discrete Mathematics 178 (1998) 229--231.


\bibitem{BBGMP98}
G. Bacs\'o, E. Boros, V. Gurvich, F. Maffray, and M. Preissmann,
On minimal imperfect graphs with circular symmetry,
Journal of Graph Theory 29:4 (1998) 209--224.

\bibitem{BG93}
E.~Boros and V.~Gurvich.
When is a circular graph minimally imperfect?
RUTCOR Research Report 22-93, Rutgers University, 1993.

\bibitem{BGZ09}
E. Boros, V. Gurvich, and  I. Zverovich.
On split and almost CIS-graphs,
Austrolasian J. of Combinatorics 43 (2009) 163--180.


\bibitem{BH09}
J. Brown, R. Hoshino,
Independence polynomials of circulants with an application to music,
Discrete Math. 309 (2009) 2292--2304.


\bibitem{BH11}
J. Brown and R. Hoshino,
Well-covered circulant graphs,
Discrete Mathematics 311 (2011) 244--251.

\bibitem{CLB81}
 {\sc D.G.~Corneil, H.~Lerchs} and {\sc L.~Stewart Burlingham}.
\newblock   Complement reducible graphs.
\newblock {\em Discrete Appl. Math.} 3 (1981) 163--174.

\bibitem{Deb56}
N. G. de Bruijn, On number systems.
Nieuw Archief voor Wiskunde 3:IV (1956) 15--17.

\bibitem{deCaen90}
D. de Caen, D. A. Gregory, I. G. Hughes, and D. L. Kreher,
Near-factors of finite groups, Ars Combinatoria 29 (1990) 53--63.

\bibitem{CGPW79}
V. Chv\'atal, R. L. Graham, A. F. Perold, and S. H. Whitesides,
Combinatorial designs related to the strong perfect graph conjecture,
Discrete Mathematics 26 (1979) 83--92.

\bibitem{CS93}
V. Chv\'atal and P.J. Slater,
A note on well-covered graphs, Ann. Discrete Math. 55 (1993) 179–182.

\bibitem{GV99}
B. Codenotti, I. Gerace, S. Vigna,
Hardness results and spectra techniques for combinatorial problems on circulant graphs,
IEEE Trans. Comput. 48 (1999) 345--351.

\bibitem{CLS81}
D.G. Corneil, H. Lerchs, L. Stewart Burlingham,
Complement reducible graphs.
Discrete Applied Mathematics 3 (1981) 163--174.


\bibitem{DengLZ04}
X. Deng, G. Li, and W. Zang,
Proof of Chv\'atal's conjecture on maximal stable sets
and maximal cliques in graphs,
J. Combin. Theory Ser. B91:2 (2004) 301--325.

\bibitem{DengLZ05}
X. Deng, G. Li, and W. Zang,
Corrigendum to: "Proof of Chv\'atal's conjecture on maximal
stable sets and maximal cliques in graphs",
[J. Combin. Theory Ser. B {\bf 91} (2) (2004) 301--325],
J. Combin. Theory Ser. B {\bf 94} (2) (2005) 352--353

\bibitem{Ful71}
D.R. Fulkerson,
Blocking and anti-blocking pairs of polyhedra,
Math. Programming 1 (1971) 168--194.

\bibitem{GGSTU93}
A.Gol'berg, V. Gurvich, A. Shapovalov, A. Temkin, and V. Udalov,
Rotational graphs without odd holes and anti-holes,
Russian Acad. Sci. Dokl. Math. 47:2 (1993) 278--284.

\bibitem{Gri84}
C. M. Grinstead, On circular critical graphs, Discrete Mathematics 51:1 (1984) 11--24.

\bibitem{Gur77}
V.Gurvich,
On repetition-free Boolean functions,
Uspechi mat. nauk (Russian Math. Surveys) 32:1  (1977) 183-184 (in Russian).

\bibitem{Gur84}
V. Gurvich,
Some properties and applications of
complete edge-chromatic graphs and hypergraphs,
Soviet math. dokl. 30:3 (1984) 803-807.


\bibitem{Gur11}
V. Gurvich,
On exact blockers and anti-blockers, $\Delta$-conjecture, and
related problems, Discrete Appl. Math. 159 (2011) 311--321.

\bibitem{GT94}
V. Gurvich and  A. Temkin,
Berge's Conjecture holds for rotational graphs;
Russian Acad. Sci. Dokl. Math. 48:2 (1994) 271--278.

\bibitem{Hoshino-thesis}
R. Hoshino, Independence Polynomials of Circulant Graphs,
Ph.D. Thesis, Dalhousie University, 2008.

\bibitem{HIK}
R. Hammack, W. Imrich, S. Klav\v zar, Handbook of product graphs. Second edition.
Discrete Mathematics and its Applications (Boca Raton). CRC Press, Boca Raton, FL, 2011.

\bibitem{Lov72a}
L.~Lov\'asz,
Normal hypergraphs and the perfect graph conjecture, Discrete Math. 2 (1972) 253--267.


\bibitem{KS06}
K. Kashiwabara and Tadashi Sakuma,
Grinstead's conjecture is true for graphs with a small clique number,
Discrete Mathematics 306 (2006) 2572--2581.


\bibitem{Muzychuk}
M. Muzychuk,
\'Ad\'am's conjecture is true in the square-free case,
J. Combin. Theory Ser. A 72 (1995) 118--134.

\bibitem{Plu70}
M. Plummer, Some covering concepts in graphs, J. Combin. Theory 8 (1970) 91--98.

\bibitem{Plu93}
M. Plummer, Well-covered graphs: a survey, Quaest. Math. 16 (1993) 253--287.

\bibitem{SS92}
R.S. Sankaranarayana, L.K. Stewart,
Complexity results for well-covered graphs, Networks 22 (1992) 247--262.

\bibitem{Sum71}
D.P. Sumner, Indecomposable graphs, Ph.D. Thesis, Univ. of
Massachuesetts, Amherst, 1971.

\bibitem{TanTar97}
D. Tankus, M. Tarsi,
The structure of well-covered graphs and the complexity of their recognition problems.
J. Comb. Theory, Ser. B 69 (1997) 230--233.

\bibitem{Uno-Mace}
T. Uno, MACE: MAximal Clique Enumerator, ver.~2.0.\\
\texttt{http://research.nii.ac.jp/$\sim$uno/code/mace.html}.

\bibitem{WZZ}
Y. Wu, W. Zang, C.-Q. Zhang,
A Characterization of Almost CIS Graphs.
SIAM J. Discrete Math. 23 (2009) 749--753.

\bibitem{Zang}
W. Zang,
Generalizations of Grillet's theorem
on maximal stable sets and maximal cliques in graphs,
Discrete Mathematics 143 (1995) 259--268.
\end{thebibliography}
\end{document}